\documentclass[a4paper,12pt,reqno]{amsart}
\usepackage[T1]{fontenc}
\usepackage[utf8]{inputenc}
\usepackage[english]{babel}
\usepackage{fullpage}
\usepackage{amsmath,amssymb}
\usepackage{ifthen}
\usepackage{arrayjobx}
\usepackage{psfrag}
\usepackage{tikz}
\usepackage[margin=0pt]{caption}
\usepackage{hyperref}
\def\HH{\mathcal H}
\def\II{\mathcal I}
\def\KK{\mathcal K}
\def\Ccea{C_\text{\rm C\'ea}}
\def\UU{\mathcal{U}}

\newcommand{\eps}{\varepsilon}

\newcommand{\normal}{\boldsymbol{n}}

\newcommand{\N}{\mathbb{N}}

\newcommand{\R}{\mathbb{R}}
\newcommand{\T}{\mathbb{T}}

\newcommand{\GG}{\mathcal{G}}

\newcommand{\MM}{\mathcal{M}}
\newcommand{\NN}{\mathcal{N}}
\newcommand{\PP}{\mathcal{P}}
\newcommand{\RR}{\mathcal{R}}
\renewcommand{\SS}{\mathcal{S}}
\newcommand{\TT}{\mathcal{T}}
\newcommand{\XX}{\mathcal{X}}

\def\Copt{C_{\rm opt}}
\def\Ctwo{C_{\rm est}}

\def\Capx{C_{\rm apx}}
\def\Cmesh{C_{\rm mesh}}
\def\Cson{C_{\rm son}}

\newcommand{\Cmark}{C_{\mathrm{mark}}}

\newcommand{\Cstab}{C_{\mathrm{stab}}}
\newcommand{\qlin}{q_{\mathrm{lin}}}
\newcommand{\qsat}{q_{\mathrm{sat}}}

\def\ksat{\kappa_{\rm sat}}
\def\kksat{\widetilde\kappa_{\rm sat}}
\DeclareMathOperator*{\diam}{diam}

\DeclareMathOperator*{\refine}{refine}
\DeclareMathOperator*{\supp}{supp}
\newcommand{\enorm}[3][]{#1|\!#1|\!#1|\,#2\,#1|\!#1|\!#1|_{#3}}

\newcommand{\norm}[3][]{#1\|#2#1\|_{#3}}
\def\d#1{\,{\rm d}#1}
\def\enorm#1{|\!|\!|#1|\!|\!|}
\def\edual#1#2{\langle\hspace*{-.6ex}\langle #1 \,,\, #2\rangle\hspace*{-.6ex}\rangle}
\def\dual#1#2{\langle #1 \,,\, #2\rangle}
\def\set#1#2{\big\{#1 \,:\, #2\big\}}
\def\reff#1#2{\stackrel{\eqref{#1}}{#2}}
\def\normal{\boldsymbol{n}}

\def\OO{\mathcal{O}}
\def\A{\mathbb{A}}
\def\H{\widetilde{H}}
\def\coarse{H}
\def\fine{h}
\def\opt{{\rm opt}}

\newcounter{statement}
\newenvironment{statement}[2][!]{%
\vskip3mm
\hrule
\hrule
\hrule
\vskip1mm
\noindent%
\refstepcounter{statement}%
\bf#2~\thestatement%
\ifthenelse{\equal{#1}{!}}{.\ }{~(#1).\ }%
\it%
}{%
\vskip1mm
\hrule
\hrule
\hrule
\vskip2mm
}
\newenvironment{theorem}[1][!]{\begin{statement}[#1]{Theorem}}{\end{statement}}

\newenvironment{proposition}[1][!]{\begin{statement}[#1]{Proposition}}{\end{statement}}
\newenvironment{corollary}[1][!]{\begin{statement}[#1]{Corollary}}{\end{statement}}
\newenvironment{remark}[1][!]{\begin{statement}[#1]{Remark}}{\end{statement}}
\newenvironment{algorithm}[1][!]{\begin{statement}[#1]{Algorithm}}{\end{statement}}

\usepackage{fancyhdr}
\advance\footskip0.4cm
\advance\textheight-0.4cm
\calclayout
\makeatletter
\def\@seccntformat#1{%
  \protect\textup{\protect\@secnumfont
    \hspace*{5mm}\ifnum\pdfstrcmp{subsection}{#1}=0 \bfseries\fi
    \csname the#1\endcsname
    \protect\@secnumpunct
  }%
}
\makeatother

\makeatletter
\def\section{\@startsection{section}{1}%
\z@{.7\linespacing\@plus\linespacing}{.5\linespacing}%
{\normalsize\scshape\bfseries\centering}}
\makeatother

\makeatletter
\renewcommand{\@secnumfont}{\bfseries}
\makeatother

\title{The saturation assumption yields\\ optimal convergence of two-level adaptive BEM}
\author{Dirk Praetorius}
\author{Michele Ruggeri}
\address{TU Wien, Institute for Analysis and Scientific Computing, Workgroup on Numerics of PDEs, Wiedner Hauptstrasse~8--10, 1040 Vienna, Austria}
\email{dirk.praetorius@asc.tuwien.ac.at}
\email{michele.ruggeri@asc.tuwien.ac.at \quad \rm(corresponding author)}
\author{Ernst P.\ Stephan}
\address{Leibniz University Hannover, Institute of Applied Mathematics, Welfengarten~1, 30167 Hannover, Germany}
\email{stephan@ifam.uni-hannover.de}
\thanks{\emph{Acknowledgements.}
The author DP acknowledges the support of the Austria Science Fund (FWF)
through the research project \emph{Optimal adaptivity for BEM and FEM-BEM coupling} (grant P27005)
and the special research program (SFB) \emph{Taming complexity in partial differential systems} (grant F65).
}

\date{\today}
\begin{document}

\begin{abstract}
We consider the convergence of adaptive BEM
for weakly-singular and hypersingular integral equations
associated with the Laplacian and the Helmholtz operator in 2D and 3D.
The local mesh-refinement is driven by some two-level error estimator.
We show that the adaptive algorithm drives the underlying error estimates to zero.
Moreover, we prove that the saturation assumption already implies linear convergence
of the error with optimal algebraic rates.
\end{abstract}
\maketitle

\def\XXX{{\color{red}\bf XXX}}
\def\HH{\mathcal{H}}
\section{Introduction}

The idea of using the difference of two approximations of different orders to obtain 
a computable estimate for the error is a well-known technique in the numerical analysis
of ordinary~\cite{pd1981} and partial differential equations~\cite[Chapter~5]{AinsworthOden00}.
Following this concept, hierarchical error estimators were among the first strategies for \textsl{a~posteriori} error estimation
of finite element and boundary element computations~\cite{bw1985,bs1993,bank1996,ms1999fembem,ms2000,hos2011,dsm2012}.
Two-level error estimators are intimately connected with, but different to, hierarchical error estimators.
A well-known disadvantage of this class of estimators is that the crucial upper error bound
(usually referred to as \emph{reliability} of the estimator) relies on a so-called \emph{saturation assumption}
(see, e.g., \eqref{eq:intro:saturation} below) and, in many situations, is even equivalent to that.
However, two-level error estimators perform strikingly well in practice;
see, e.g., \cite{bek1993,mms1997,msw1998,eh2006,effp2009,efgp2013}.
The purpose of the present work is to shed some light on this empirical observation.

Let us illustrate the concept of two-level error estimation with a concrete example.
For instance, we consider the weakly-singular integral equation
\begin{equation} \label{eq:intro:weaksing}
(Vu)(x) := \int_\Gamma G(x-y) \, u(y) \d{y} = f(x)
\quad \text{for (almost) all } x \in \Gamma := \partial \Omega
\end{equation}
associated with the Laplace operator $-\Delta$ in 2D.
Here, $\Gamma$ is the boundary of a bounded Lipschitz domain $\Omega \subset \R^2$,
$f \in H^{1/2}(\Gamma)$ is a given right-hand side,
$G(\cdot)$ is the fundamental solution of $-\Delta$,
and $u \in H^{-1/2}(\Gamma)$ is the sought integral density;
see Section~\ref{section:weaksing:setting} for more details on the precise functional analytic setting.
It is well-known that, for certain right-hand sides $f$, this integral equation is equivalent
to the homogeneous Laplace equation posed on the domain $\Omega$
supplemented with some inhomogeneous Dirichlet boundary conditions on $\Gamma$.
Without loss of generality, we may assume that $\diam(\Omega) < 1$.
Then, the variational formulation of~\eqref{eq:intro:weaksing} seeks $u \in \HH := H^{-1/2}(\Gamma)$ such that
\begin{equation}\label{eq:intro:continuous}
 \dual{Vu}{v} := \int_\Gamma (Vu)(x) \, v(x) \d{x} 
 = \int_\Gamma f(x) v(x) \d{x} =: \dual{f}{v}
 \quad \text{for all } v \in \HH.
\end{equation}
We note that the left-hand side defines a scalar product on the energy space $\HH$,
while the right-hand side defines a linear and continuous functional on $\HH$.
We denote by $\enorm{v}^2 := \dual{Vv}{v}$ the operator-induced energy norm.
Given a partition $\TT_h$ of $\Gamma$ into line segments
and the associated finite-dimensional subspace
$\XX_h := \set{w_h \in L^2(\Gamma)\!}{\!w_h|_T \text{ constant for all } T \in \TT_h} \subset \HH$,
a possible conforming Galerkin discretization reads
\begin{equation} \label{eq:intro:galerkin}
\dual{Vu_h}{v_h} = \dual{f}{v_h}
\quad \text{ for all } v_h \in \XX_h.
\end{equation}
The Lax--Milgram lemma proves existence and uniqueness of both the continuous solution $u \in \HH$ 
of~\eqref{eq:intro:continuous} and the Galerkin approximation $u_h \in \XX_h$ of~\eqref{eq:intro:galerkin}.

In this 2D setting, the two-level error estimator reads as follows:
Let $\XX_{h/2}$ be the finite-dimensional space associated with the partition $\TT_{h/2}$ of $\Gamma$
obtained by bisecting all $T \in \TT_h$.
Let $u_{h/2} \in \XX_{h/2}$ be the corresponding Galerkin solution.
For each $T \in \TT_h$, let $\widehat\varphi_T \in \XX_{h/2}$ be the fine-mesh Haar function
(which takes the values $+1$ and $-1$)
such that $\supp(\widehat\varphi_T) = T$ and $\int_T \widehat\varphi_T \d{x} = 0$.
Then, one can show that, for some $C>0$, it holds that
\begin{equation} \label{eq:intro:aposteriori}
 C^{-1} \, \enorm{u_{h/2} - u_h}^2 
 \le \sum_{T \in \TT_h} \tau_h(T)^2
 := \sum_{T \in \TT_h} \frac{|\dual{f-Vu_h}{\widehat\varphi_T}|^2}{\enorm{\widehat\varphi_T}^2}
 \le C \, \enorm{u_{h/2} - u_h}^2.
\end{equation}
The proof of the estimates in~\eqref{eq:intro:aposteriori} is based on the
direct space decomposition $\XX_{h/2} = \XX_h \oplus \bigoplus_{T \in \TT_h} {\rm span}\{\widehat\varphi_T\}$,
which is even stable in $\HH$.
In explicit terms,
the difference $\enorm{u_{h/2} - u_h}$ between the two Galerkin solutions is equivalent
to the coarse-mesh residual $f - Vu_h$ tested by the additional basis functions of $\XX_{h/2}$.
This provides a computable measure for the error improvement in the Galerkin orthogonality
\begin{equation}\label{eq:intro:orthogonality}
 \enorm{u-u_{h/2}}^2 = \enorm{u - u_h}^2 - \enorm{u_{h/2} - u_h}^2.
\end{equation}
Under the saturation assumption
\begin{equation}\label{eq:intro:saturation}
 \enorm{u-u_{h/2}} \le q \, \enorm{u - u_h}
 \quad \text{for some constant } 0 < q < 1,
\end{equation}
combining~\eqref{eq:intro:aposteriori} and~\eqref{eq:intro:orthogonality},
we thus obtain the {\sl a~posteriori} error estimate
\begin{equation*}
C^{-1} \, \sum_{T \in \TT_h} \tau_h(T)^2
\le \enorm{u-u_h}^2
\le \frac{C}{1-q^2} \, \sum_{T \in \TT_h} \tau_h(T)^2.
\end{equation*}
Subject to the saturation assumption~\eqref{eq:intro:saturation},
one can thus use the local contributions $\tau_h(T)$ as error indicators to steer an adaptive mesh-refinement algorithm
of the usual form
\begin{equation*}
 \boxed{~\rm SOLVE~}
 \quad \longrightarrow \quad
 \boxed{~\rm ESTIMATE~}
 \quad \longrightarrow \quad
 \boxed{~\rm MARK~}
 \quad \longrightarrow \quad
 \boxed{~\rm REFINE~}
\end{equation*}
which empirically leads to striking numerical results, even though one cannot expect that the saturation assumption~\eqref{eq:intro:saturation} holds in general~\cite{bek1996}.

In the context of FEM (instead of BEM), it follows from the analysis in~\cite{msv2008} that the usual adaptive algorithm
drives the two-level error estimator to zero.
Moreover, in this context, the sum of two-level error estimator and data oscillations is locally equivalent to
the usual residual error estimator.
Therefore, it follows from~\cite{ks2011,axioms} that one even gets optimal algebraic convergence rates
for the sum of estimator and oscillations
(even without the saturation assumption). 

In the context of BEM, the work~\cite{ffmps2014} shows that
the usual adaptive algorithm drives the two-level error estimator to zero.
However, the analysis is restricted to the Laplace equation and relies on the weighted-residual error estimator
from~\cite{cs1996,cc1997,cms2001} as an auxiliary tool
in combination with non-trivial inverse estimates from~\cite{invest}. 

In the present work, we simplify the argument from~\cite{ffmps2014}
and prove that the adaptive algorithm leads to estimator convergence for (quite general) integral equations.
In addition, under an appropriate variant of the saturation assumption~\eqref{eq:intro:saturation}
(see~\eqref{eq:saturation} in Section~\ref{sec:saturation} below),
we prove that the adaptive algorithm guarantees even linear convergence with optimal algebraic convergence rates
for the energy error.
This analysis is developed in an abstract framework in the spirit of~\cite{axioms}
and covers weakly-singular as well as hypersingular integral equations associated with the Laplacian
as well as the Helmholtz operator in 2D and 3D.

The remainder of this work is organized as follows:
Section~\ref{section:twolevel} proposes
an abstract setting for two-level {\sl a~posteriori} error estimation
(Sections~\ref{section:twolevel:abstract}--\ref{section:estimator}),
formulates the adaptive algorithm (Algorithm~\ref{algorithm}),
and proves linear convergence (Theorem~\ref{theorem:linear})
as well as optimal convergence rates (Theorem~\ref{theorem:optimal}).
In Section~\ref{section:weaksing}, we show that our abstract framework covers weakly-singular integral equations
(with energy space $\H^{-1/2}(\Gamma)$) for the 2D and 3D Laplace and Helmholtz equation
(Theorem~\ref{proposition:weaksing:twolevel}).
Moreover, for the same setting and without requiring any saturation assumption,
we prove plain estimator convergence (Theorem~\ref{theorem:weaksing:plain}).
Finally, Section~\ref{section:hypsing} shows that the abstract framework also covers hypersingular integral equations
(with energy space $\H^{1/2}(\Gamma)$) for the 2D and 3D Laplace and Helmholtz equation
(Theorem~\ref{proposition:hypsing:twolevel}).
Also for this case, without requiring any saturation assumption,
we prove plain estimator convergence (Theorem~\ref{theorem:hypsing:plain}).

\section{An abstract analysis of two-level adaptivity}
\label{section:twolevel}

\subsection{Abstract problem and its discretization}
\label{section:twolevel:abstract}
Let $\HH$ be a Hilbert space with scalar product $\edual\cdot\cdot$ and corresponding norm $\enorm\cdot$.
Let $\KK: \HH \to \HH'$ be a compact operator, where $\HH'$ is the dual space of $\HH$. We denote the duality brackets on $\HH' \times \HH$ by $\dual\cdot\cdot$. Given $F \in \HH'$, we suppose that the variational formulation
\begin{equation}\label{eq:abstract:weakform}
 b(u,v) := \edual{u}{v} + \dual{\KK u}{v} = F(v)
 \quad\text{for all } v \in \HH
\end{equation}
admits a unique solution $u \in \HH$. 

Let $\T$ denote the set of all admissible triangulations.
For each $\TT_\coarse \in \T$, let $\XX_\coarse \subset \HH$ be the associated conforming subspace.
The corresponding Galerkin discretization of~\eqref{eq:abstract:weakform} reads:
Find $u_\coarse \in \XX_\coarse$ such that
\begin{equation}\label{eq:abstract:galerkin}
 b(u_\coarse, v_\coarse) = F(v_\coarse)
 \quad\text{for all } v_\coarse \in \XX_\coarse.
\end{equation}
We assume that there exists $\gamma_0 > 0$ with
\begin{equation}\label{eq1:infsup}
\inf_{w_\coarse \in \XX_\coarse \setminus \{0\}} 
\sup_{v_\coarse \in \XX_\coarse \setminus \{0\}} 
\frac{|b(w_\coarse, v_\coarse)|}{\enorm{w_\coarse} \, \enorm{v_\coarse}}
\ge \gamma_0
\quad \text{for all } \TT_\coarse \in \T.
\end{equation}
In particular,~\eqref{eq:abstract:galerkin} admits a unique solution $u_\coarse \in \XX_\coarse$, and there holds the C\'ea lemma
\begin{equation}\label{eq:abstract:cea}
 \enorm{u - u_\coarse} \le \Ccea \, \min_{v_\coarse \in \XX_\coarse} \enorm{u - v_\coarse},
\end{equation}
where
\begin{equation*}
 \Ccea := 1 + \frac{M}{\gamma_0}
 \quad \text{and} \quad 
 \displaystyle M := \sup_{v, w \in \HH \setminus \{0\}} \frac{| b(w,v)|}{\enorm{w} \, \enorm{v}} < \infty.
\end{equation*}
Finally, we assume that all $u \in \HH$ can be approximated by discrete functions, i.e.,
\begin{equation}\label{eq:approximation-property}
 \inf_{\TT_\coarse \in \T} \min_{v_\coarse \in \XX_\coarse} \enorm{v - v_\coarse} = 0
 \quad \text{ for all $v \in \HH$.}
\end{equation}

\begin{remark}
If $\KK = 0$, then~\eqref{eq1:infsup}--\eqref{eq:abstract:cea} hold with $\gamma_0 = 1 = \Ccea$ and equality.
\end{remark}

\subsection{Mesh-refinement}
Given a fixed mesh-refinement algorithm,
for $\TT_\coarse \in \T$ and marked elements $\MM_\coarse \subseteq \TT_\coarse$,
we denote by $\TT_\fine := \refine(\TT_\coarse, \MM_\coarse) \in \T$ the coarsest refinement of $\TT_\coarse$
such that all marked elements are refined, i.e., $\MM_\coarse \subseteq \TT_\coarse \setminus \TT_\fine$.
We write $\TT_\fine \in \refine(\TT_\coarse)$, if $\TT_\fine \in \T$ is obtained from $\TT_\coarse \in \T$ by finitely many steps of refinement.
Moreover, we denote by $\widehat\TT_\coarse := \refine(\TT_\coarse,\TT_\coarse) \in \T$
the uniform refinement of $\TT_\coarse$.
Throughout, we assume that refinement leads to nested discrete spaces,
i.e., 
\begin{equation}\label{eq:nestedness}
 \XX_\coarse \subseteq \XX_\fine
 \quad \text{for all } \TT_\coarse \in \T \text{ and } \TT_\fine \in \refine(\TT_\coarse).
\end{equation}%
We assume that $\T = \refine(\TT_0)$ for a fixed initial triangulation $\TT_0$,
i.e., admissible triangulations are refinements of $\TT_0$.
Moreover, we suppose that, for all $\TT_\coarse, \TT_\fine \in \T$
satisfying $\TT_\fine = \refine(\TT_\coarse, \MM_\coarse)$ for some $\MM_\coarse \subseteq \TT_\coarse$,
it holds that
\begin{enumerate}
\renewcommand{\theenumi}{{M\arabic{enumi}}}
\setcounter{enumi}{0}
\bf
\item\label{assumption:sons}\rm
$T = \bigcup\set{T' \in \TT_\fine}{T'\subseteq T}$ and $\#\set{T'\in\TT_\fine}{T' \subseteq T}\le\Cson$
\quad
for all $T \in \TT_\coarse$,
\end{enumerate}
where $\Cson \ge 2$ is a uniform constant.

In addition to the foregoing general assumptions, we require the following additional properties:
For all $\TT,\TT' \in \T$, there exists a coarsest common refinement $\TT \oplus \TT' \in \refine(\TT) \cap \refine(\TT')$ such that
\begin{enumerate}
\renewcommand{\theenumi}{{M\arabic{enumi}}}
\setcounter{enumi}{1}
\bf
\item\label{assumption:overlay}\rm
$\#(\TT \oplus \TT') \le \#\TT + \#\TT' - \#\TT_0$.
\end{enumerate}
Moreover, there exists $\Cmesh > 0$, which depends only on $\TT_0$, such that for all sequences of triangulations such that $\TT_{\ell + 1} = \refine(\TT_\ell, \MM_\ell)$ with arbitrary $\MM_\ell \subseteq \TT_\ell$, it holds that
\begin{enumerate}
\renewcommand{\theenumi}{{M\arabic{enumi}}}
\setcounter{enumi}{2}
\bf
\item\label{assumption:closure}\rm
$ \#\TT_\ell - \#\TT_0
 \le \Cmesh \, \sum_{j = 0}^{\ell - 1} \#\MM_j
 \quad \text{for all } \ell \in \N_0$.
\end{enumerate}

Assumption~\eqref{assumption:sons} specifies that one-level refinement does only lead to a bounded number of son elements and that parental elements are the union of their children. It is proved in~\cite{gss2014} for refinement by newest vertex bisection (NVB) and holds with $\Cson = 4$ in 2D.
Assumption~\eqref{assumption:overlay} is called \emph{overlay estimate}.
It is first found in~\cite{stevenson2007,ckns2008} for NVB.
Assumption~\eqref{assumption:closure} is called \emph{closure estimate}.
It is first proved in~\cite{bdd2004,stevenson} for NVB assuming an admissibility condition on the initial mesh $\TT_0$.
For NVB in 2D, the admissibility assumption on $\TT_0$ from~\cite{bdd2004,stevenson}
is proved to be unnecessary in~\cite{kpp}.
For a 1D bisection algorithm which ensures~\eqref{assumption:sons}--\eqref{assumption:closure} together with boundedness of the local mesh-ratio
\begin{equation*}
 |T|/|T'| \le \gamma < \infty
 \quad \text{for all } \TT \in \T \text{ and } T, T' \in \TT \text{ with } T \cap T' \neq\emptyset,
\end{equation*}
we refer to~\cite{cmam2013}.
Moreover,~\eqref{assumption:sons}--\eqref{assumption:closure} also holds for red-refinement with first-order hanging nodes for meshes consisting of triangles/simplices or rectangles/cuboids;
see~\cite{bn2010}.
Finally, for results on mesh-refinement in the frame of isogeometric analysis, we refer to~\cite{mp2015,bgmp2016,ghp2017}.

\begin{remark}\label{remark:approximation}
{\rm(a)} By definition of $\refine(\cdot)$ and~\eqref{eq:nestedness}, assumption~\eqref{eq:approximation-property} is equivalent to the assumption that uniform mesh-refinement leads to convergence
of the best approximation error to zero.

{\rm(b)} The uniform discrete inf-sup condition~\eqref{eq1:infsup} is compatible
with the assumptions on the mesh-refinement in the following sense:
If~\eqref{eq:abstract:weakform} admits a unique solution
and uniform mesh-refinement leads to convergence of the best approximation error to zero,
then~\eqref{eq1:infsup} follows from the nestedness of the discrete spaces, if $\TT_0$ is sufficiently fine;
see, e.g.,~\cite[Proposition~1]{bhp2017}.
\end{remark}

\subsection{Two-level error estimation}
\label{section:estimator}

For $\TT_\coarse \in \T$, let $\II_\coarse$ be the index set
corresponding to the degrees of freedom of the space $\XX_\coarse$
associated with $\TT_\coarse$.
Similarly, let $\widehat\II_\coarse$ be the index set
corresponding to the degrees of freedom of the space $\widehat\XX_\coarse$
associated with the uniform refinement $\widehat\TT_\coarse = \refine(\TT_\coarse,\TT_\coarse)$ of $\TT_\coarse$.
By nestedness of the discrete spaces, it holds that $\II_\coarse \subseteq \widehat\II_\coarse$.
For each $z \in \widehat\II_\coarse \setminus \II_\coarse$,
we denote by $\tau_\coarse(z) \geq 0$ the associated error indicator.
For $T \in \TT_\coarse$, let $\widehat\II_\coarse(T) \subseteq \widehat\II_\coarse$ be the corresponding degrees of freedom with respect to the fine mesh $\TT_\fine$, i.e., we have the (in general non-disjoint) union
\begin{equation*}
\widehat\II_\coarse \setminus \II_\coarse = \bigcup_{T \in \TT_{\coarse}} \widehat\II_\coarse(T).
\end{equation*}
We define the two-level element error indicators by
\begin{subequations} \label{eq:twolevel}
\begin{equation}
 \tau_\coarse(T)^2 := \sum_{z \in \widehat\II_\coarse(T) \setminus \II_\coarse} \tau_\coarse(z)^2
\quad
\text{for all } T \in \TT_\coarse.
\end{equation}
The resulting two-level element-based error estimator $\tau_\coarse$ is given by
\begin{equation}
 \tau_\coarse := \tau_\coarse(\TT_\coarse),
 \quad \text{where} \quad
 \tau_\coarse(\UU_\coarse) := \bigg(\sum_{T \in \UU_\coarse} \tau_\coarse(T)^2 \bigg)^{1/2}
 \quad \text{for all } \UU_\coarse \subseteq \TT_\coarse.
\end{equation}
\end{subequations}
For the element-based estimator~\eqref{eq:twolevel}, we make the following assumption:
There exists $\Ctwo > 0$ such that, for all $\TT_\coarse \in \T$, $\MM_\coarse \subseteq \TT_\coarse$,
and $\TT_\fine := \refine(\TT_\coarse, \MM_\coarse)$, it holds that

\begin{enumerate}
\renewcommand{\theenumi}{{E\arabic{enumi}}}
\setcounter{enumi}{0}
\bf
\item\label{assumption:reliable}\rm
$\displaystyle\Ctwo^{-1} \, \tau_\coarse(\MM_\coarse)
\le \enorm{u_\fine - u_\coarse} 
\le \Ctwo \, \tau_\coarse(\TT_\coarse \setminus \TT_\fine)$,
\end{enumerate}
where we recall that $\MM_\coarse \subseteq \TT_\coarse \setminus \TT_\fine$. Finally, we suppose that there exists $\Cstab > 0$ such that, for all $\TT_\coarse \in \T$ and all $\TT_\fine \in \refine(\TT_\coarse)$, it holds that
\begin{enumerate}
\renewcommand{\theenumi}{{E\arabic{enumi}}}
\setcounter{enumi}{1}
\bf
\item\label{assumption:stable}\rm
$| \tau_\fine(\TT_\coarse \cap \TT_\fine) - \tau_\coarse(\TT_\coarse \cap \TT_\fine) |
\le \Cstab \, \enorm{u_\fine - u_\coarse}$.
\end{enumerate}

The lower estimate in~\eqref{assumption:reliable} is usually named \emph{discrete efficiency} and was first exploited in~\cite{doerfler1996,mns2000}. Moreover, together with the saturation assumption (see Section~\ref{sec:saturation} below), the upper bound in~\eqref{assumption:reliable} will yield \emph{(discrete) reliability} of the two-level error estimator. 
Assumption~\eqref{assumption:stable} is usually referred to as \emph{stability} of the error estimator, i.e., the estimator on non-refined elements depends Lipschitz continuously on the discrete solutions~\cite{axioms}. 

\subsection{Adaptive algorithm}

We consider the following standard adaptive algorithm. Under the assumptions~\eqref{assumption:sons}--\eqref{assumption:closure} and~\eqref{assumption:reliable}--\eqref{assumption:stable}, we will show that an appropriate saturation assumption~\eqref{eq:saturation} yields linear convergence with optimal algebraic rates.

\begin{algorithm}\label{algorithm}
{\bfseries Input:} Initial mesh $\TT_0$, adaptivity parameters $0 < \theta \le 1$ and $1 \le \Cmark \le \infty$.
\\
{\bfseries Loop:} For all $\ell = 0, 1, 2, \dots$, iterate the following steps~{\rm(i)--(iv)}:
\begin{itemize}
\item[\rm(i)] 
Compute the discrete solution $u_\ell \in \XX_\ell$.
\item[\rm(ii)] 
Compute the two-level indicators $\tau_{\ell}(T)$ for all $T \in \TT_\ell$.
\item[\rm(iii)] 
Determine $\MM_\ell \subseteq \TT_\ell$ of almost minimal cardinality
(i.e., minimal up to the multiplicative factor $\Cmark$),
which satisfies the D\"orfler marking criterion~\cite{doerfler1996}
\begin{equation}\label{eq:doerfler}
 \theta \, \tau_\ell^2  
 \le \tau_{\ell}(\MM_\ell)^2.
\end{equation}
\item[\rm(iv)] Define $\TT_{\ell+1} := \refine(\TT_\ell,\MM_\ell)$, increase the counter $\ell \mapsto \ell+1$, and goto {\rm(i)}.
\end{itemize}
{\bfseries Output:} Meshes $\TT_\ell$, corresponding solutions $u_\ell$, and error estimators $\tau_\ell$ for all $\ell \in \N_0$.
\end{algorithm}

\begin{remark}
{\rm(a)}
If $\MM_\ell^{\mathrm{min}} \subseteq \TT_\ell$ denotes the (in general non-unique) set of minimal cardinality
which satisfies~\eqref{eq:doerfler}, almost minimal cardinality of $\MM_\ell$ means that
\begin{equation} \label{eq:doerfler:minimal}
\# \MM_\ell \leq \Cmark \, \# \MM_\ell^{\mathrm{min}}.
\end{equation}
For an algorithm with linear complexity determining $\MM_\ell$
satisfying \eqref{eq:doerfler} and \eqref{eq:doerfler:minimal} with $\Cmark=2$,
we refer to~\cite{stevenson2007}.
Moreover, we refer to the recent work~\cite{pp2019} 
for an algorithm with linear cost $\OO(\#\TT_\ell)$ and $\Cmark = 1$.

{\rm(b)}
The choice $\Cmark = \infty$ means that $\MM_\ell \subseteq \TT_\ell$ satisfies the marking criterion~\eqref{eq:doerfler}, but may be (too) large. In particular, uniform mesh-refinement $\MM_\ell = \TT_\ell$ is allowed for all $0 < \theta \le 1$.

{\rm(c)}
If $\Cmark < \infty$, small $0 < \theta \ll 1$ generically leads to few marked elements and thus highly adapted meshes, while only $\theta = 1$ leads to uniform mesh-refinement $\MM_\ell = \TT_\ell$.

\end{remark}

\subsection{Plain convergence}

In concrete situations, Algorithm~\ref{algorithm} guarantees estimator convergence $\tau_\ell \to 0$ as $\ell \to \infty$.
For instance, we refer
to Theorem~\ref{theorem:weaksing:plain} in Section~\ref{section:weaksing} for weakly-singular integral equations
and to Theorem~\ref{theorem:hypsing:plain} in Section~\ref{section:hypsing} for hypersingular integral equations.
Under an appropriate saturation assumption, however, one can even prove much stronger convergence results.

\subsection{Saturation assumption} \label{sec:saturation}

Let $0 < \ksat \le \qsat < 1$.
We suppose the following assumption along the sequence of adaptive meshes generated by Algorithm~\ref{algorithm}:
For all $\ell \in \N_0$ and all refinements $\TT_\fine \in \refine(\TT_\ell)$, the mesh
$\TT_\coarse := \refine(\TT_\ell, \TT_\ell \setminus \TT_\fine)$ satisfies that
\begin{enumerate}
\renewcommand{\theenumi}{{S}}
\bf
\item\label{eq:saturation}\rm
if $\enorm{u - u_\fine} \le \ksat \, \enorm{u - u_\ell}$, then
$\enorm{u - u_\coarse} \le \qsat \, \enorm{u - u_\ell}$.
\end{enumerate}
In explicit terms, the saturation assumption~\eqref{eq:saturation} states that,
if $\TT_\fine \in \refine(\TT_\ell)$ leads to a sufficient improvement of the error,
then the one-level refinement of $\TT_\coarse$ towards $\TT_\fine$ already provides a uniform improvement
of the error.

\begin{remark}
{\rm(a)} In the literature (see, e.g.,~\cite{mms1997,msw1998,ms2000,hms2001,heuer2002,eh2006,fp2008,effp2009,efgp2013,affkp2015}), the saturation assumption is usually formulated with respect to the uniform refinement, and it is assumed that
 \begin{equation}\label{eq:old:saturation}
  \enorm{u - \widehat u_\ell} \le \qsat \, \enorm{u - u_\ell}
  \quad \text{for all } \ell \in \N_0
 \end{equation}
 holds with a uniform constant $0 < \qsat < 1$ along the sequence of adaptively generated meshes. In contrast to~\eqref{eq:saturation}, the previous saturation assumption~\eqref{eq:old:saturation} states that uniform one-level refinement provides a uniform improvement of the error.

{\rm(b)} In fact,~\eqref{eq:saturation} implies~\eqref{eq:old:saturation}. To see this, note that each sufficiently fine mesh $\TT_\fine \in \refine(\TT_\ell)$ yields that $\enorm{u - u_\fine} \le \ksat \, \enorm{u - u_\ell}$ and $\TT_\ell \cap \TT_\fine = \emptyset$. Hence, it follows that $\TT_\coarse = \refine(\TT_\ell, \TT_\ell \setminus \TT_\fine) = \widehat\TT_\ell$, and~\eqref{eq:saturation} thus implies that $\enorm{u - \widehat u_\ell} \le \qsat \, \enorm{u - u_\ell}$.

{\rm(c)} We will see below that the weaker saturation assumption~\eqref{eq:old:saturation} yields linear convergence (Theorem~\ref{theorem:linear}). However, our proof of optimal convergence rates (Theorem~\ref{theorem:optimal}) relies on the stronger saturation assumption~\eqref{eq:saturation}.

{\rm(d)} It is well-known that counterexamples show that even the weaker saturation assumption~\eqref{eq:old:saturation} cannot hold in general~\cite{bek1996,dn2002}. However, as soon as the numerical scheme has reached its asymptotic behavior, then at least~\eqref{eq:old:saturation} is valid; see the discussion in~\cite[Section~5.2]{fp2008}.
\end{remark}

\subsection{Linear convergence}

Our first observation is that the saturation assumption~\eqref{eq:saturation}
(or at least its weaker form~\eqref{eq:old:saturation})
already yields linear convergence of the error.

\begin{theorem}[linear convergence]\label{theorem:linear}
Suppose assumption~\eqref{assumption:reliable}.
Consider the output of Algorithm~\ref{algorithm} for fixed parameters $0 < \theta \le 1$ and $1 \le \Cmark \le \infty$.
We suppose that the saturation assumption~\eqref{eq:saturation} holds (or at least the weaker form~\eqref{eq:old:saturation}).
Then, there exist $0 < \qlin < 1$ and $\ell_0 \in \N_0$ such that
Algorithm~\ref{algorithm} yields linear convergence of the energy error
\begin{equation}\label{eq1:linear}
 \enorm{u - u_{\ell + 1}} \le \qlin \, \enorm{u - u_\ell}
 \quad \text{for all } \ell \ge \ell_0.
\end{equation}
The index $\ell_0$ depends only on $\KK$, $u$, as well as the sequence $(u_\ell)_{\ell \in \N_0}$ of discrete solutions.
The constant $\qlin$ depends only on $\ell_0$, $\Ctwo$, $\qsat$, and $\theta$.
\end{theorem}

\begin{proof}
The proof structurally follows the ideas developed in the early works~\cite{doerfler1996,mns2000}.
Note that $\MM_\ell \subseteq \TT_\ell \setminus \TT_{\ell + 1}$.
With the lower bound of~\eqref{assumption:reliable} (used for $\TT_\fine = \TT_{\ell+1}$), the marking strategy~\eqref{eq:doerfler} implies that
\begin{equation}\label{eq1:proof:linear}
 \Ctwo^2 \, \enorm{u_{\ell + 1} - u_{\ell}}^2
 \reff{assumption:reliable}\ge \tau_\ell(\MM_\ell)^2
 \reff{eq:doerfler}\ge \theta \, \tau_\ell^2.
\end{equation}
Together with the saturation assumption~\eqref{eq:old:saturation}, the triangle inequality proves that
\begin{equation*}
 \enorm{u - u_\ell}
 \le \enorm{u - \widehat u_\ell} + \enorm{\widehat u_\ell - u_\ell}
 \le \qsat \, \enorm{u - u_\ell} + \enorm{\widehat u_\ell - u_\ell}.
\end{equation*}
With the upper bound of~\eqref{assumption:reliable} (used for $\TT_\fine = \widehat\TT_\ell$), it follows that
\begin{equation} \label{eq:abstract:convergence1}
 (1 - \qsat)^2 \, \enorm{u - u_{\ell}}^2 
 \le \enorm{\widehat u_\ell - u_{\ell}}^2
 \reff{assumption:reliable}\le \Ctwo^2 \tau_\ell^2 
 \reff{eq1:proof:linear}\le \frac{\Ctwo^4}{\theta} \, \enorm{u_{\ell+1} - u_\ell}^2.
\end{equation}
According to~\cite[Lemma~4.2]{msv2008}, there exists $u_\infty \in \XX_\infty := \overline{\bigcup_{\ell=0}^\infty \XX_\ell}$ such that
\begin{equation*}
 \enorm{u_\infty - u _\ell} \to 0
 \quad \text{as } \ell \to \infty.
\end{equation*}
Hence, the right-hand side of the last estimate tends to zero as $\ell \to \infty$. This implies that $u = u_\infty$ and allows to apply~\cite[Lemma~18]{bhp2017}: For all $0 < \eps < 1$, there exists an index $\ell_0 \in \N_0$ such that
\begin{equation}\label{eq:abstract:pythagoras}
 \enorm{u - u_{\ell+1}}^2 + \enorm{u_{\ell+1} - u_\ell}^2 \le \frac{1}{1 - \eps} \, \enorm{u - u_\ell}^2
 \quad \text{for all } \ell \ge \ell_0.
\end{equation}
With this Pythagoras-type estimate, we obtain that
\begin{equation*}
 \enorm{u - u_{\ell + 1}}^2 
 \reff{eq:abstract:pythagoras}{\le}
 \frac{1}{1 - \eps} \, \enorm{u - u_{\ell}}^2 - \enorm{u_{\ell + 1} - u_{\ell}}^2
 \reff{eq:abstract:convergence1}{\le}
 \Big(\frac{1}{1 - \eps}  - \frac{\theta \, (1-\qsat)^2}{\Ctwo^4} \Big) \, \enorm{u - u_{\ell}}^2.
\end{equation*}
Choosing $\eps > 0$ sufficiently small and $\ell_0 \in \N_0$ accordingly, we get $0 \le \qlin < 1$ such that
\begin{equation*}
 \enorm{u - u_{\ell + 1}}^2 \le \qlin^2 \, \enorm{u - u_{\ell}}^2
 \quad \text{for all } \ell \ge \ell_0.
\end{equation*}
This concludes the proof of~\eqref{eq1:linear}. 
\end{proof}

\begin{corollary}
If $\KK = 0$, then Theorem~\ref{theorem:linear} holds with $\ell_0 = 0$, and the constant $\qlin$ depends only on $\Ctwo$, $\qsat$, and $\theta$.
\end{corollary}

\begin{proof}
Note that $\edual{u - u_{\ell+1}}{v_{\ell+1}} = 0$ for all $v_{\ell+1} \in \XX_{\ell+1}$ and, in particular, for $v_{\ell+1} = u_{\ell+1} - u_\ell$.
Instead of~\eqref{eq:abstract:pythagoras}, we can thus employ the Pythagoras theorem
\begin{equation*}
 \enorm{u - u_{\ell+1}}^2 + \enorm{u_{\ell+1} - u_\ell}^2 = \enorm{u - u_\ell}^2
 \quad \text{for all } \ell \in \N_0.
\end{equation*}
This even simplifies the proof.
\end{proof}

\begin{remark}
{\rm(a)} Note that the weaker saturation assumption~\eqref{eq:old:saturation} is stated with respect to the uniform refinement $\widehat\TT_\ell$ of $\TT_\ell$. Then, the (structurally identical) implication~\eqref{eq1:linear} of Theorem~\ref{theorem:linear} states the uniform improvement of the energy error with respect to the adaptively refined mesh $\TT_{\ell + 1}$. The interpretation is that Algorithm~\ref{algorithm} guarantees sufficient enrichment of $\TT_\ell$ in each step of the adaptive loop (if such an enrichment exists).

{\rm(b)} If $\KK = 0$, then the saturation assumption~\eqref{eq:old:saturation} is equivalent to reliability $\enorm{u - u_\ell} \lesssim \enorm{\widehat u_\ell - u_\ell} \simeq \tau_\ell$ of the two-level error estimator; see also, e.g.,~\cite{fp2008}. 
In particular,~\eqref{eq1:linear} cannot be proved without the saturation assumption.

{\rm(c)} We note that linear convergence is independent of the assumptions~\eqref{assumption:sons}--\eqref{assumption:closure} on the mesh-refinement and, moreover, does not rely on stability~\eqref{assumption:stable} of the two-level error estimator.
\end{remark}

\subsection{Optimal algebraic rates}

For $N \in \N$, define the (finite) set of meshes
\begin{equation*}
 \T(N) := \set{\TT_\coarse \in \T}{\#\TT_\coarse - \#\TT_0\le N}.
\end{equation*}
For $s > 0$, define the approximation constant
\begin{equation}\label{eq:As}
 \norm{u}{\A_s}
 := \sup_{N \in \N_0} \min_{\TT_\opt \in \T(N)} \min_{v_\opt \in \XX_\opt} (N+1)^s \, \enorm{u - v_\opt} \in [0,\infty].
\end{equation}
By definition, $\norm{u}{\A_s} < \infty$ implies that the best approximation error decays with algebraic rate $s > 0$ along a sequence of optimal meshes, i.e., $\min_{v_\opt \in \XX_\opt} \enorm{u - v_\opt} = \OO\big((\#\TT_\opt)^{-s}\big)$, if $\TT_\opt$ attains the minimum in~\eqref{eq:As}.

The following theorem states that Algorithm~\ref{algorithm} guarantees that each possible algebraic rate $s > 0$ in $\enorm{u-u_\ell} = \OO\big((\#\TT_\ell)^{-s}\big)$ will, in fact, be realized along the sequence of adaptively generated meshes. 
Unlike the proof of linear convergence (Theorem~\ref{theorem:linear}), the proof relies on the stronger saturation assumption~\eqref{eq:saturation}.

\begin{theorem}[rate optimality]\label{theorem:optimal}
Suppose assumptions~\eqref{assumption:sons}--\eqref{assumption:closure} and~\eqref{assumption:reliable}--\eqref{assumption:stable}.
Consider the output of Algorithm~\ref{algorithm} for fixed parameters $0 < \theta \le 1$ and $1 \le \Cmark < \infty$.
We suppose that the saturation assumption~\eqref{eq:saturation} holds.
Then, there exists $0 < \theta_\opt \le 1$ such that, if $0 < \theta < \theta_\opt$, it holds that
\begin{equation}\label{eq:optimal}
 \forall \, s > 0 \, \exists \, \Copt > 0 \quad
 \Copt^{-1} \, \norm{u}{\A_s}
 \le \sup_{\ell \in \N_0} (\#\TT_\ell)^s \enorm{u - u_\ell}
 \le \Copt \, \norm{u}{\A_s}.
\end{equation}
The constant $\Copt$ in~\eqref{eq:optimal} depends only on $\theta$, $s$, $\ksat$, $\qsat$, $\TT_0$, and the index $\ell_0$ from Theorem~\ref{theorem:linear}.
\end{theorem}

\def\Crel{C_{\rm rel}}
\def\Cgal{C_{\rm gal}}
\def\Cq{C(q)}
\begin{proof}
For all $s>0$ such that $\norm{u}{\A_s}=0$, \eqref{eq:optimal} is trivially satisfied.
Therefore, we assume that $\norm{u}{\A_s} \in (0,\infty].$
The reminder of the proof is split into five steps.

{\bf Step~1 (discrete reliability).}
Let $\TT_\fine \in \refine(\TT_\ell)$ and $\TT_\coarse := \refine(\TT_\ell, \TT_\ell \setminus \TT_\fine)$.
In this step, we prove that $\enorm{u - u_\coarse} \le \qsat \, \enorm{u - u_\ell}$ implies that
\begin{equation}\label{eq:reliability}
 \enorm{u - u_\ell}^2 
 \le \frac{\enorm{u_\coarse - u_\ell}^2}{(1-\qsat)^2}
 \le \frac{\Ctwo^2}{(1-\qsat)^2} \, \tau_\ell(\TT_\ell \setminus \TT_\fine)^2.
\end{equation}
To see this, note that 
$\TT_\ell \setminus \TT_\coarse = \TT_\ell \setminus \TT_\fine$. Hence,~\eqref{assumption:reliable} leads to
\begin{equation*}
 (1-\qsat)^2 \, \enorm{u - u_\ell}^2
 \le \enorm{u_\coarse - u_\ell}^2
 \reff{assumption:reliable}\le \Ctwo^2 \, \tau_\ell( \TT_\ell \setminus \TT_\coarse )^2
 = \Ctwo^2 \, \tau_\ell( \TT_\ell \setminus \TT_\fine )^2.
\end{equation*}

{\bf Step~2 (optimality of D\"orfler marking).}
Define $\Crel := \Ctwo \, (\Ccea + 1) / (1 - \qsat)$.
For given $0 < \theta < \theta_{\rm opt} := (1+\Cstab^2 \Crel^2)^{-1}$, choose $\delta > 0$ and $0 < \kksat \le \ksat$ sufficiently small such that
\begin{equation}\label{eq:theta}
 \theta 
 \le \frac{1-(1+\delta^{-1}) \, \Ctwo^2 \, \Crel^2 \, \kksat^{\,2}}{1+(1+\delta) \, \Cstab^2 \Crel^2}
 < \frac{1}{1 + \Cstab^2 \Crel^2 } = \theta_\opt.
\end{equation}
Let $\ell  \in \N_0$ and $\TT_\fine \in \refine(\TT_\ell)$ such that
$\enorm{u - u_\fine} \le \kksat \, \enorm{u - u_\ell}$.
We prove that 
\begin{equation}\label{eq:optimal_doerfler}
 \theta \, \tau_\ell^2 \le \tau_\ell(\TT_\ell \setminus \TT_\fine)^2.
\end{equation}
To this end, note that $\kksat \le \ksat$ and the saturation assumption~\eqref{eq:saturation} allows us to employ~\eqref{eq:reliability} from Step~1. Since
$\XX_\ell \subseteq \XX_\fine$, the C\'ea lemma~\eqref{eq:abstract:cea} yields that
\begin{equation}\label{eq2502}
 \enorm{u_\fine \!-\! u_\ell}
 \le \enorm{u \!-\! u_\fine} + \enorm{u \!-\! u_\ell}
 \reff{eq:abstract:cea}\le (\Ccea + 1) \, \enorm{u \!-\! u_\ell}
 \reff{eq:reliability}\le \Crel \, \tau_\ell( \TT_\ell \setminus \TT_\fine ).
\end{equation}
The same argument proves that
\begin{equation} \label{eq:2806}
\enorm{\widehat u_\fine - u_\fine} \reff{eq:abstract:cea}\le (\Ccea + 1) \, \enorm{u - u_\fine}.
\end{equation}
Therefore, we obtain that
\begin{equation}\label{eq:monotonicity}
\begin{split}
\tau_\fine 
& \reff{assumption:reliable}\le \Ctwo \enorm{\widehat u_\fine - u_\fine}
\reff{eq:2806}\le \Ctwo \, (\Ccea + 1) \, \enorm{u - u_\fine}
\le  \Ctwo \, (\Ccea + 1) \, \kksat \, \enorm{u - u_\ell}
 \\
 & \reff{eq:reliability}\le \frac{\Ctwo \, (\Ccea + 1)}{1-\qsat} \, \kksat \, \enorm{\widehat u_\ell - u_\ell}
 = \Crel \, \kksat \, \enorm{\widehat u_\ell - u_\ell}
 \reff{assumption:reliable}\le \Ctwo \, \Crel \, \kksat \, \tau_\ell. 
 \end{split}
\end{equation}
Together with stability~\eqref{assumption:stable}, the Young inequality yields that
\begin{equation*}
\begin{split}
&\tau_\ell^2 
 \stackrel{\phantom{\eqref{assumption:stable}}}{=}
 \tau_\ell( \TT_\ell \setminus \TT_\fine )^2 + \tau_\ell(\TT_\ell \cap \TT_\fine)^2
 \\& \quad
 \reff{assumption:stable}\le
 \tau_\ell( \TT_\ell \setminus \TT_\fine )^2
  + (1+\delta) \, \Cstab^2 \, \enorm{u_\fine - u_\ell}^2
  + (1+\delta^{-1}) \, \tau_\fine(\TT_\ell \cap \TT_\fine)^2
 \\& \quad
 \reff{eq2502}\le \big( 1 + (1+\delta) \, \Cstab^2 \Crel^2  \big) \, \tau_\ell( \TT_\ell \setminus \TT_\fine )^2 
 + (1+\delta^{-1}) \, \tau_\fine^2
 \\& \quad
 \reff{eq:monotonicity}\le \big( 1 + (1+\delta) \, \Cstab^2 \Crel^2 \big) \, \tau_\ell( \TT_\ell \setminus \TT_\fine )^2 
 + (1+\delta^{-1}) \, \Ctwo^2 \Crel^2 \, \kksat^{\,2} \, \tau_\ell^2.
\end{split}
\end{equation*}
Rearranging this estimate and using~\eqref{eq:theta}, we prove~\eqref{eq:optimal_doerfler}. We note that $\kksat$ depends only on $\ksat$, $\qsat$, $\Ctwo$, $\Ccea$, $\theta$, and $\TT_0$.

{\bf Step~3 (comparison lemma).} Let $\ell \in \N_0$.
We show that there exists a set $\RR_\ell \subseteq \TT_\ell$ such that
\begin{equation}\label{eq:comparison}
 \#\RR_\ell \le \Capx^{1/s} \, \norm{u}{\A_s}^{1/s} \enorm{u-u_\ell}^{-1/s}
 \quad \text{and} \quad
 \theta \, \tau_\ell^2 \le \tau_\ell(\RR_\ell)^2,
\end{equation}
where $\Capx := \Ccea/\kksat > 0$.
To see this, note that, without loss of generality, we may assume that
$\norm{u}{\A_s} < \infty$ and $\enorm{u-u_\ell} > 0$ (otherwise, choose $\RR_\ell = \TT_\ell$).
Define
\begin{equation}\label{eq:def:eps}
0 < \eps := \kksat \, \Ccea^{-1} \, \enorm{u-u_\ell}
\reff{eq:abstract:cea}\le \kksat \min_{v_0 \in \XX_0} \enorm{u-v_0}
< \min_{v_0 \in \XX_0} \enorm{u-v_0} \reff{eq:As}\le \norm{u}{\A_s}
\end{equation}
and note that, by construction, it holds that $\norm{u}{\A_s}^{1/s} \, \eps^{-1/s} > 1$.
Choose $N \in \N$ minimal such that $1 \leq N < \norm{u}{\A_s}^{1/s} \, \eps^{-1/s} \le N + 1$.
Choose $\TT_\eps \in \T(N)$ and $v_\eps \in \XX_\eps$ such that $\enorm{u-v_\eps} = \min_{\TT_\opt \in \T(N)} \min_{v_\opt \in \XX_\opt} \enorm{u-v_\opt}$.
Finally, define $\RR_\ell := \TT_\ell \setminus \TT_\fine$, where $\TT_\fine := \TT_\eps \oplus \TT_\ell$ is the overlay from~\eqref{assumption:overlay}.
Altogether, we obtain that
\begin{equation*}
\begin{split}
\# \RR_\ell &= \#(\TT_\ell \setminus \TT_\fine) \reff{assumption:sons}\le \#\TT_\fine - \#\TT_\ell 
 \reff{assumption:overlay}\le \#\TT_\eps - \#\TT_0 \le N
 \\
 & < \norm{u}{\A_s}^{1/s} \, \eps^{-1/s}
 = \kksat^{-1/s} \, \Ccea^{1/s} \, \norm{u}{\A_s}^{1/s} \, \enorm{u-u_\ell}^{-1/s},
\end{split}
\end{equation*}
which is the first inequality in~\eqref{eq:comparison}.
Moreover, $\TT_\fine \in \refine(\TT_\eps)$ and the choice of $N$ prove that
\begin{equation*}
\begin{split}
\Ccea^{-1} \, \enorm{u-u_\fine}
& \reff{eq:abstract:cea}\le \enorm{u-v_\eps} =\!\!\! \min_{\TT_\opt \in \T(N)} \min_{v_\opt \in \XX_\opt} \enorm{u-v_\opt}
 \\&
 \reff{eq:As}\le (N+1)^{-s} \norm{u}{\A_s}
 \le  \eps
 \reff{eq:def:eps}= \kksat \, \Ccea^{-1} \, \enorm{u-u_\ell}
\end{split}
\end{equation*}
and thus $\enorm{u - u_\fine} \le \kksat \, \enorm{u - u_\ell}$.
Hence, we can use~\eqref{eq:optimal_doerfler} from Step~2, which yields
the second inequality in \eqref{eq:comparison} and thus concludes the proof.

{\bf Step~4 (adaptivity guarantees optimal rates).} Since Algorithm~\ref{algorithm} chooses $\MM_\ell \subseteq \TT_\ell$ with essentially minimal cardinality, it follows that 
\begin{equation}\label{eq2:comparison}
\# \MM_\ell \reff{eq:doerfler:minimal}\le \Cmark \# \RR_\ell \le \Cmark \Capx^{1/s} \, \norm{u}{\A_s}^{1/s} \enorm{u-u_\ell}^{-1/s}
\quad\text{ for all $\ell \in \N_0$.}
\end{equation}
We recall from~\cite[Lemma~22]{bhp2017} that, for $\ell > 0$, it holds that
\begin{equation}\label{eq:bhp2017}
 \#\TT_\ell \!-\! \#\TT_0 
 < \#\TT_\ell \!-\! \#\TT_0 \!+\! 1
 \le \#\TT_\ell
 \le \#\TT_0 \, (\#\TT_\ell \!-\! \#\TT_0 \!+\! 1)
 \le 2 \, \#\TT_0 \, (\#\TT_\ell \!-\! \#\TT_0).
\end{equation}
Let $\ell_0 \in \N_0$ be the index from Theorem~\ref{theorem:linear}.
Define 
\begin{equation}\label{eq:const}
 C := \max_{j = 0,\dots,\ell_0} \frac{\#\MM_j}{\#\MM_{\ell_0}} \ge 1.
\end{equation}
Recall from Theorem~\ref{theorem:linear} that $\enorm{u-u_j}^{-1/s} \le \qlin^{(j-\ell)/s} \, \enorm{u-u_\ell}^{-1/s}$ for all $\ell_0 \le j \le \ell$.
For $\ell > \ell_0$, the geometric series thus proves that
\begin{equation*}
\begin{split}
 &\#\TT_\ell \simeq \#\TT_\ell - \#\TT_0
 \reff{assumption:closure}\lesssim \sum_{j = 0}^{\ell - 1} \#\MM_j
 = \sum_{j = 0}^{\ell_0-1} \#\MM_j + \sum_{j = \ell_0}^{\ell-1} \#\MM_j
 \\& \quad
 \reff{eq:const}\le (\ell_0 C + 1) \, \sum_{j = \ell_0}^{\ell-1} \#\MM_j
 \reff{eq2:comparison}\lesssim \norm{u}{\A_s}^{1/s} \sum_{j=\ell_0}^{\ell - 1} \enorm{u-u_j}^{-1/s}
 \reff{eq1:linear}\le \frac{1}{1-\qlin^{1/s}} \, \norm{u}{\A_s}^{1/s} \enorm{u-u_\ell}^{-1/s}.
\end{split}
\end{equation*}
For $\ell = 0$, the estimate $\#\TT_0 \lesssim 1 \le \norm{u}{\A_s}^{1/s} \enorm{u-u_0}^{-1/s}$ is trivial by definition of $\norm{u}{\A_s}$ in~\eqref{eq:As}.
For $0 \le \ell < \ell_0$, assumption~\eqref{assumption:sons} yields that $\#\TT_\ell \le \Cson^\ell \#\TT_0 \lesssim \norm{u}{\A_s}^{1/s} \enorm{u-u_0}^{-1/s}
\le \Ccea^{1/s} \, \norm{u}{\A_s}^{1/s} \enorm{u-u_\ell}^{-1/s}$. Hence, we have proved that
$\sup_{\ell \in \N_0} (\#\TT_\ell)^s \,  \enorm{u-u_\ell} \lesssim \norm{u}{\A_s}$, which is the upper estimate in~\eqref{eq:optimal}.

{\bf Step~5 (adaptivity constraints optimal rates).} To prove the lower estimate in~\eqref{eq:optimal}, we may assume that the upper bound is finite, i.e., $\sup_{\ell \in \N_0} (\#\TT_\ell)^s \,  \enorm{u-u_\ell} < \infty$. Let $N \in \N_0$ and choose the maximal $\ell \in \N_0$ such that $\#\TT_\ell - \#\TT_0 \le N$, i.e., $\TT_\ell \in \T(N)$. Besides~\eqref{eq:bhp2017}, note that~\eqref{assumption:sons} implies that $\#\TT_\ell \le \# \TT_{\ell+1} \le \Cson \, \#\TT_\ell$.
In particular, it follows that
\begin{equation*}
 N + 1 < \#\TT_{\ell+1} - \#\TT_0 + 1 \reff{eq:bhp2017}\simeq \#\TT_{\ell+1} \simeq \#\TT_\ell.
\end{equation*}
Hence, we obtain that
\begin{equation*}
 \min_{\TT_\opt \in \T(N)} (N+1)^s \enorm{u-u_\opt} \lesssim (\#\TT_\ell)^s \, \enorm{u-u_\ell} \le \sup_{\ell \in \N_0} (\#\TT_\ell)^s \, \enorm{u-u_\ell}.
\end{equation*}
Taking the supremum over all $N \in \N_0$, we prove the lower estimate in~\eqref{eq:optimal}.
\end{proof}

\section{Weakly-singular integral equation}
\label{section:weaksing}

\subsection{Functional analytic framework for Laplace BEM}
\label{section:weaksing:setting}

Let $\Gamma \subseteq \partial\Omega$ be a relatively open and connected part of the boundary $\partial\Omega$ of a bounded Lipschitz domain $\Omega \subset \R^d$ with $d = 2,3$.
For the ease of presentation, we assume that $\partial\Omega$ is polygonal.
For $d = 3$, we additionally assume that $\Gamma$ is a Lipschitz dissection~\cite[pp.~99]{mclean}.
For $d = 2$, we additionally assume that $\diam(\Omega) < 1$, which can always be achieved by scaling.
With the fundamental solution of the Laplacian, i.e., $G(z) := -\frac{1}{2\pi} \, \log|z|$ for $d = 2$ resp.\ $G(z) := \frac{1}{4\pi} \, |z|^{-1}$ for $d=3$, we consider the weakly-singular integral equation 
\begin{equation}\label{eq:Vu=f}
 (V u) (x) :=  \int_\Gamma u(y) \, G(x-y) \d{y} = f(x)
 \quad \text{for } x \in \Gamma,
\end{equation}
where $f \in H^{1/2}(\Gamma)$ is some given right-hand side and $u \in \H^{-1/2}(\Gamma)$ is the sought integral density.
Here, $H^{1/2}(\Gamma) := \set{v|_\Gamma}{v \in H^1(\Omega)}$ is the trace space of $H^1(\Omega)$ and $\H^{-1/2}(\Gamma) := H^{1/2}(\Gamma)'$ is its dual space with respect to the (extended) $L^2(\Gamma)$-scalar product.
 
The variational formulation~\eqref{eq:abstract:weakform} of~\eqref{eq:Vu=f} reads:
Find $u \in \HH := \H^{-1/2}(\Gamma)$ such that
\begin{equation}\label{eq:weaksing}
 b(u,v) := \edual{u}{v} := \dual{Vu}{v} = \dual{f}{v}
 \quad \text{for all } v \in \HH,
\end{equation}
where $\dual{v}{w} := \int_\Gamma vw \d{x}$ denotes the $L^2(\Gamma)$-scalar product. It is well-known that the weakly-singular integral operator $V: \H^{-1/2}(\Gamma) \to H^{1/2}(\Gamma)$ is a symmetric and elliptic isomorphism. 
Therefore, $\edual\cdot\cdot$ is a scalar product, and the induced norm $\enorm{v}^2 := \edual{v}{v}$ is an equivalent norm on $\H^{-1/2}(\Gamma)$. In particular, the Lax--Milgram theorem proves the existence and uniqueness of the solution $u \in \H^{-1/2}(\Gamma)$ of~\eqref{eq:weaksing}.

\subsection{Functional analytic framework for Helmholtz BEM}
\label{section:weaksing:setting:helmholtz}

We employ the notation from Section~\ref{section:weaksing:setting}.
For a wavenumber $\kappa > 0$, let 
\begin{equation*}
 G_\kappa(z) := \frac{i}{4} \, H_0^{(1)}(\kappa|z|)
 \quad \text{for $d = 2$} \quad \text{resp.} \quad 
 G_\kappa(z) := \frac{e^{i\kappa \, |z|}}{4\pi \, |z|}
 \quad \text{for $d = 3$},
\end{equation*}
where $H_0^{(1)}$ is the first-kind Hankel function of order zero. 
We consider the weakly-singular integral equation 
\begin{equation}\label{eq:Vu=f:helmholtz}
 (V_\kappa u) (x) :=  \int_\Gamma u(y) \, G_\kappa(x-y) \d{y} = f(x)
 \quad \text{for } x \in \Gamma,
\end{equation}
where $f \in H^{1/2}(\Gamma)$ is some given right-hand side and $u \in \H^{-1/2}(\Gamma)$ is the sought integral density.
The variational formulation~\eqref{eq:abstract:weakform} of~\eqref{eq:Vu=f:helmholtz} reads:
Find $u \in \HH := \H^{-1/2}(\Gamma)$ such that
\begin{equation}\label{eq:weaksing:helmholtz}
 b(u,v) := \dual{V_\kappa u}{v} = \dual{f}{v}
 \quad \text{for all } v \in \HH.
\end{equation}
It is known that the single-layer operator $V_\kappa: \H^{-1/2}(\Gamma) \to H^{1/2}(\Gamma)$ is an isomorphism, if and only if $\kappa^2$ is not an eigenvalue of the interior Dirichlet problem for the Laplace operator; see, e.g.,~\cite[Theorem~3.9.1]{sauter-schwab}. We suppose throughout that this is the case, i.e., $V_\kappa$ is an isomorphism and~\eqref{eq:weaksing:helmholtz} thus admits a unique solution $u \in \H^{-1/2}(\Gamma)$.

\begin{subequations}\label{eq:shape-regular}
\subsection{Mesh-refinement in~2D}
\label{section:weaksing:mesh2d}
For $d = 2$, a mesh $\TT_\coarse$ is a finite partition of $\Gamma$ into compact affine line segments. We employ the bisection algorithm from~\cite{cmam2013}. We assume that all meshes $\TT_\coarse \in \T = \refine(\TT_0)$ are obtained by applying this mesh-refinement strategy to a given initial mesh $\TT_0$. As already mentioned in Section~\ref{section:twolevel}, this guarantees~\eqref{assumption:sons}--\eqref{assumption:closure}. Moreover, there holds uniform $\gamma$-shape regularity in the sense of
\begin{equation}\label{eq:1d:shape-regular}
 \sup_{\TT_\coarse \in \T}\max\set{{\rm diam}(T) / {\rm diam}(T')}{T, T' \in \TT_\coarse \text{ with } T \cap T' \neq \emptyset}
 \le \gamma < \infty,
\end{equation}
where $\gamma > 0$ depends only on $\TT_0$. 
We recall that $\widehat\TT_\coarse = \refine(\TT_\coarse; \TT_\coarse)$ denotes the uniform refinement, where all elements have been bisected once.

\begin{figure}[t]
 \centering
 \psfrag{T0}{}
 \psfrag{T1}{}
 \psfrag{T2}{}
 \psfrag{T3}{}
 \psfrag{T4}{}
 \psfrag{T12}{}
 \psfrag{T34}{}
 \includegraphics[width=35mm]{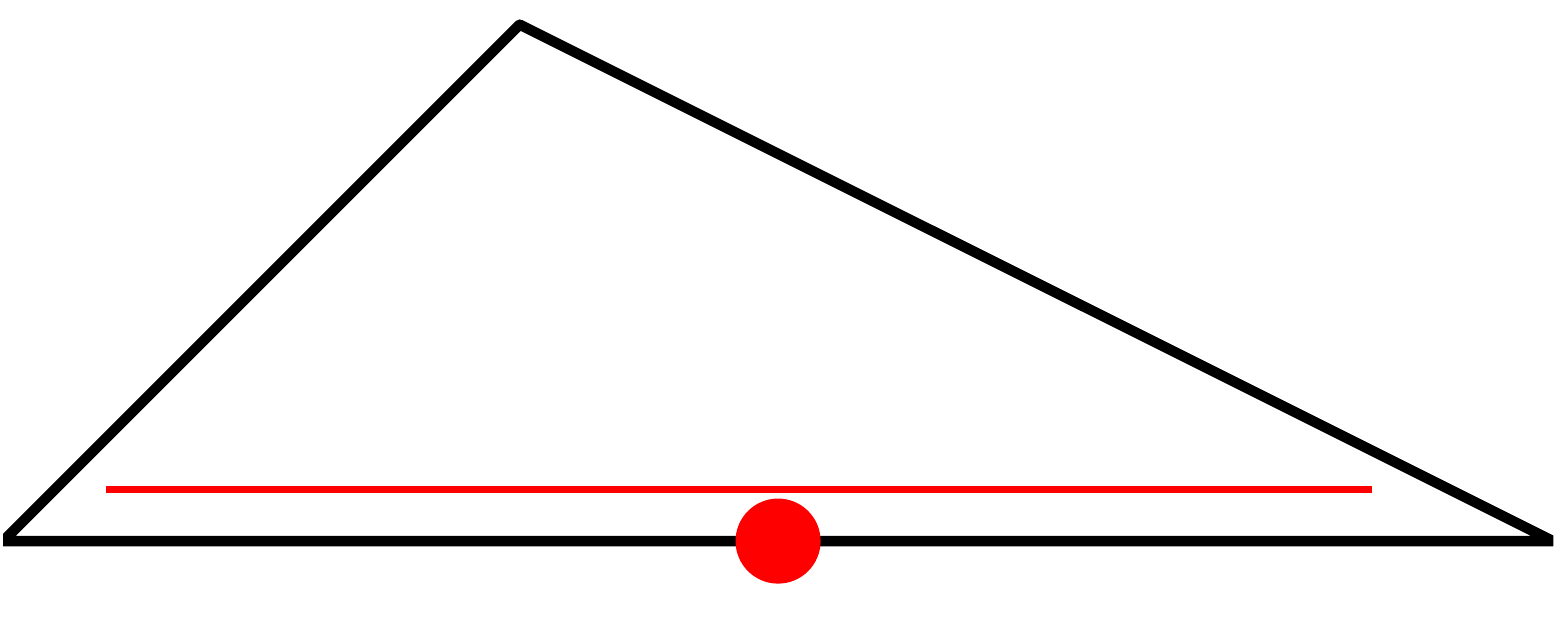} \quad
 \includegraphics[width=35mm]{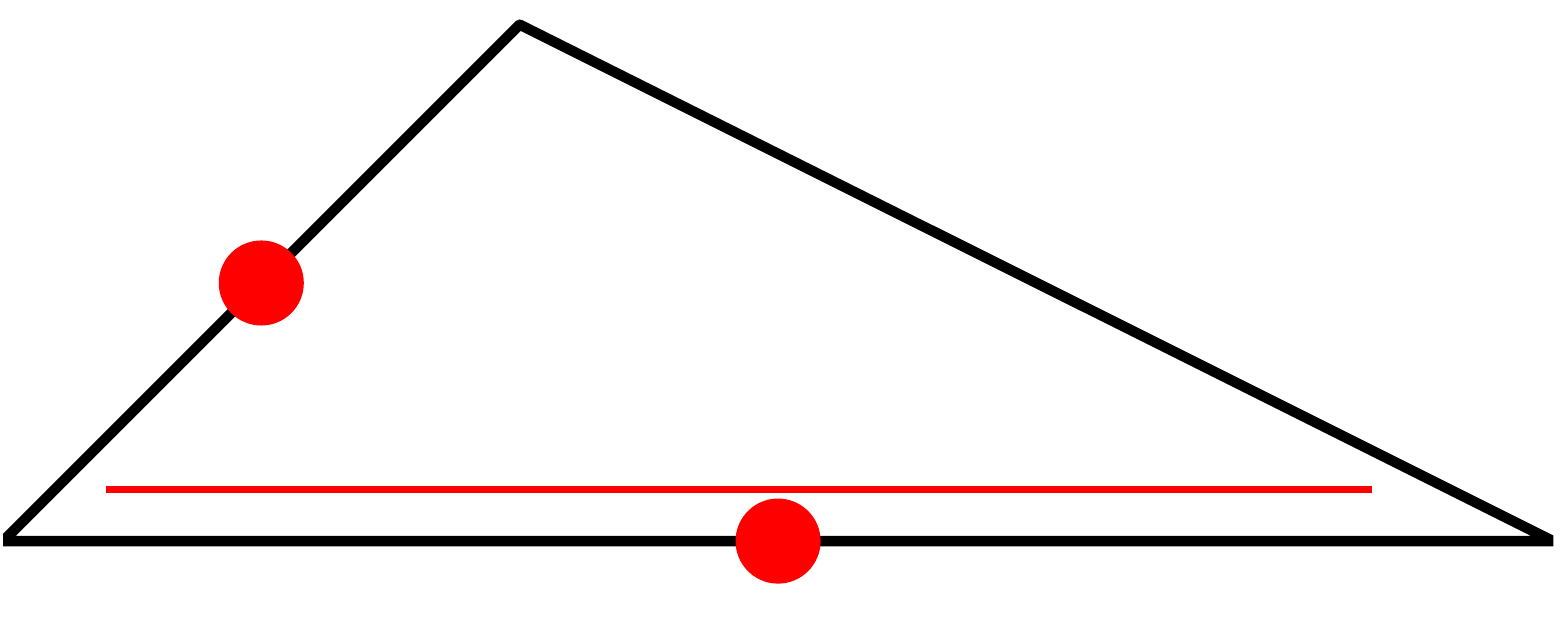} \quad
 \includegraphics[width=35mm]{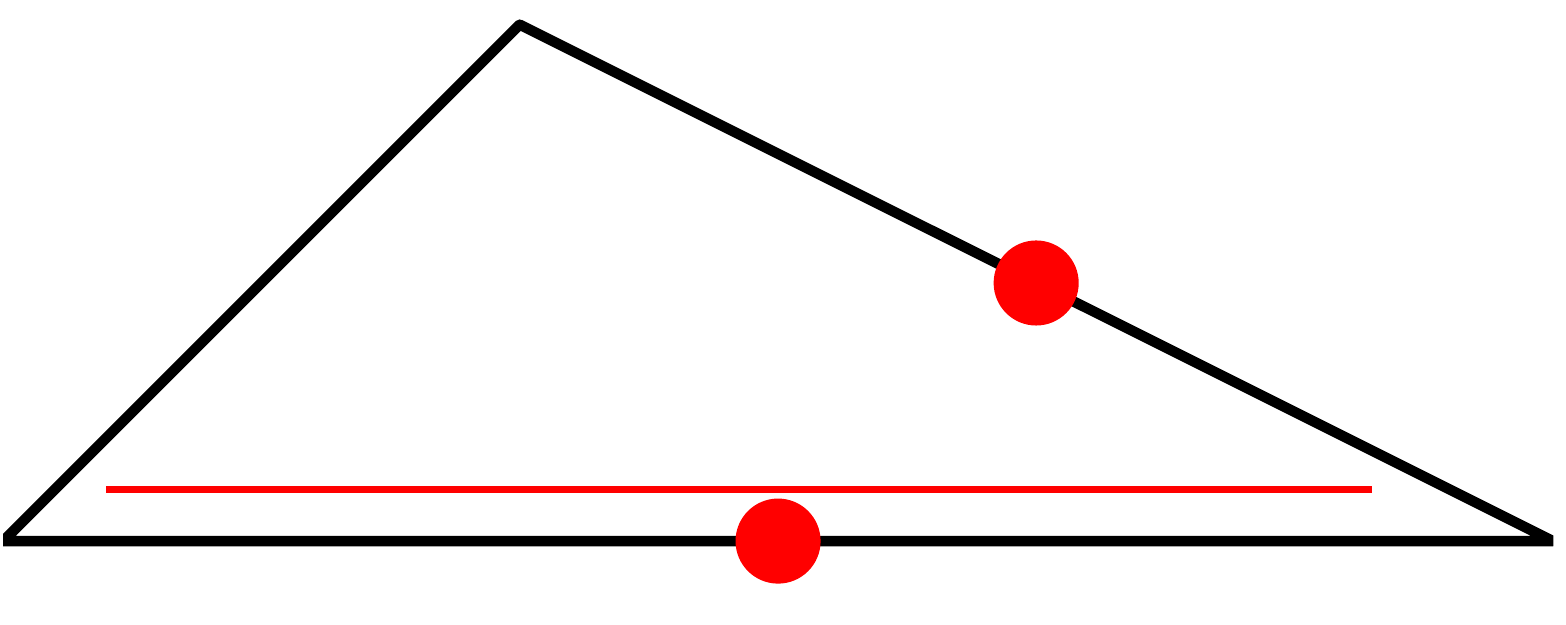} \quad
 \includegraphics[width=35mm]{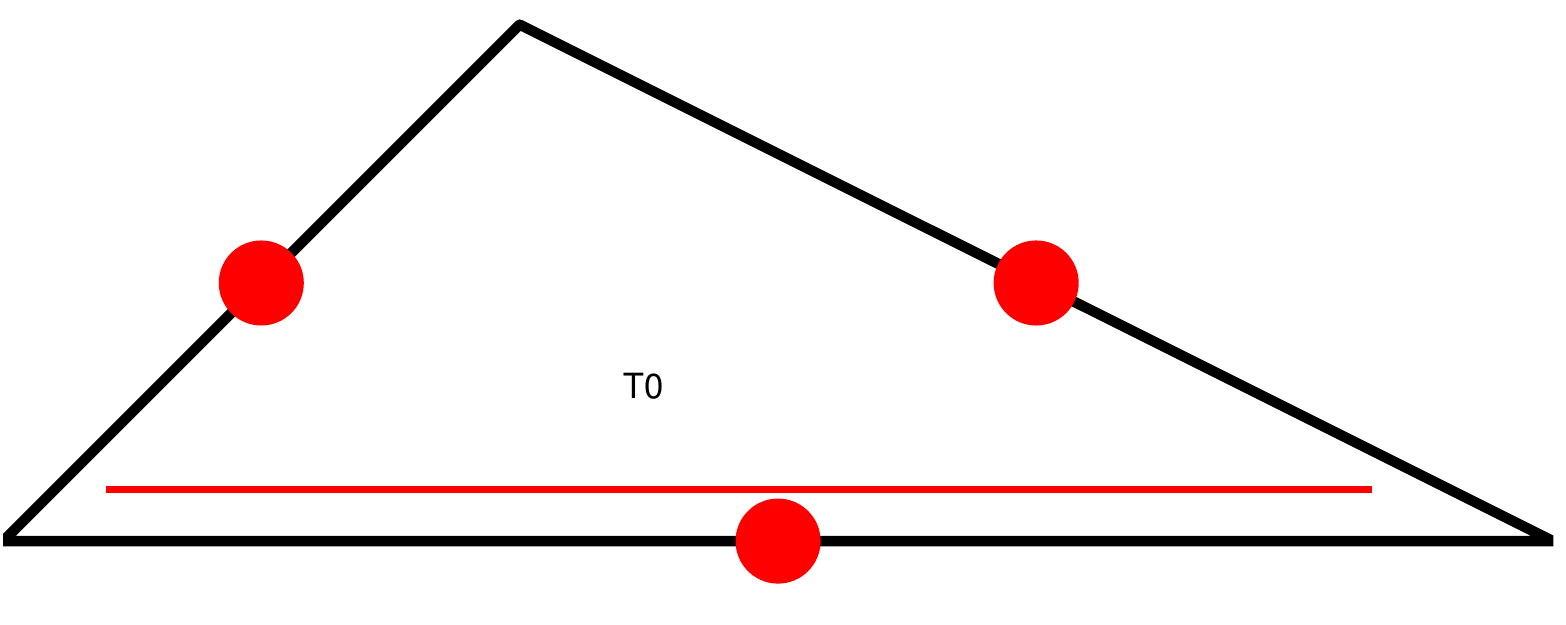} \\
 \includegraphics[width=35mm]{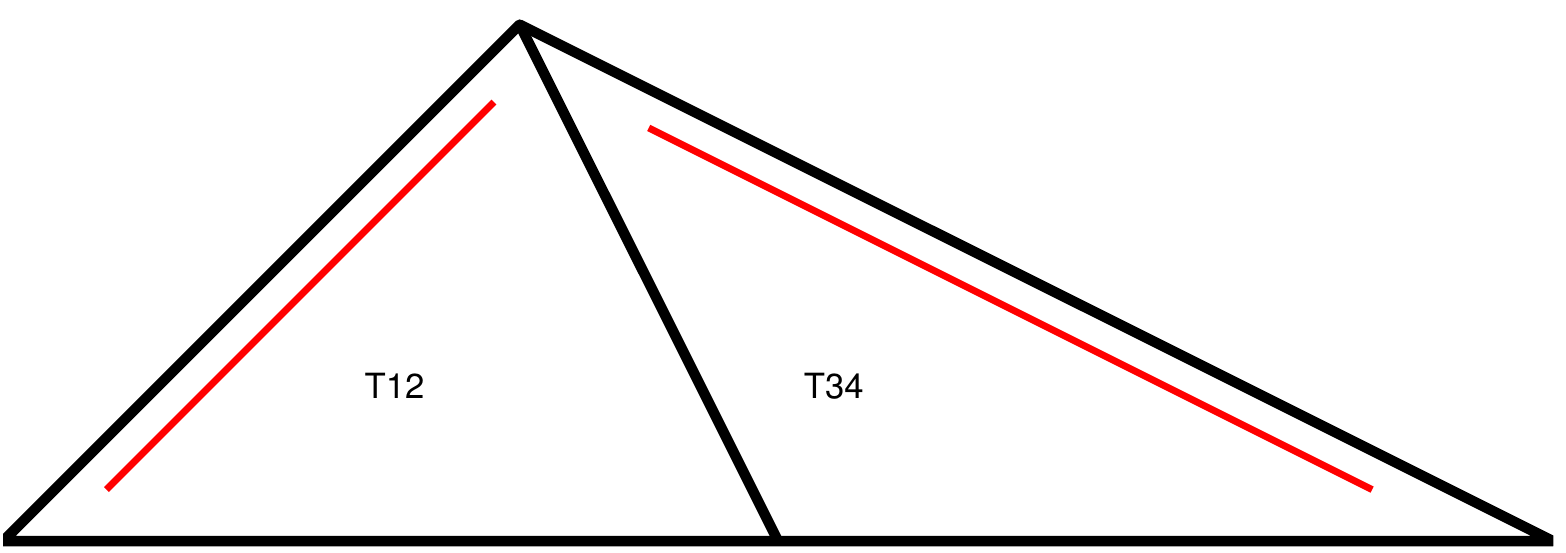} \quad
 \includegraphics[width=35mm]{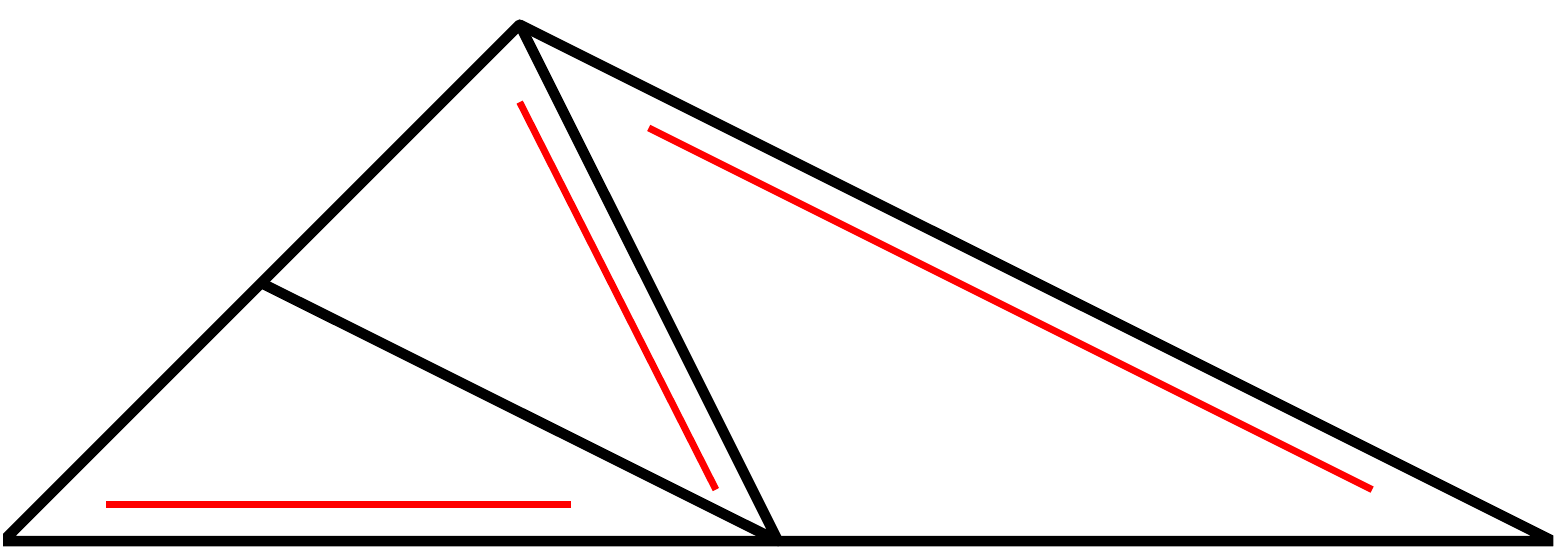}\quad
 \includegraphics[width=35mm]{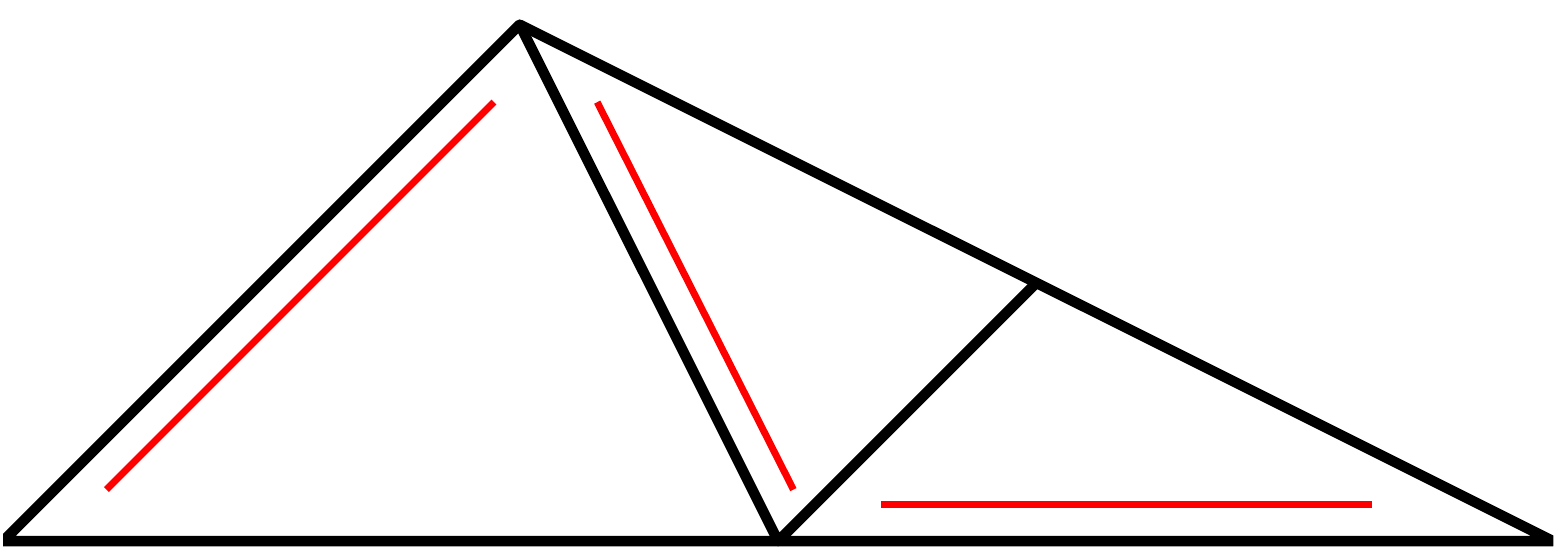}\quad
 \includegraphics[width=35mm]{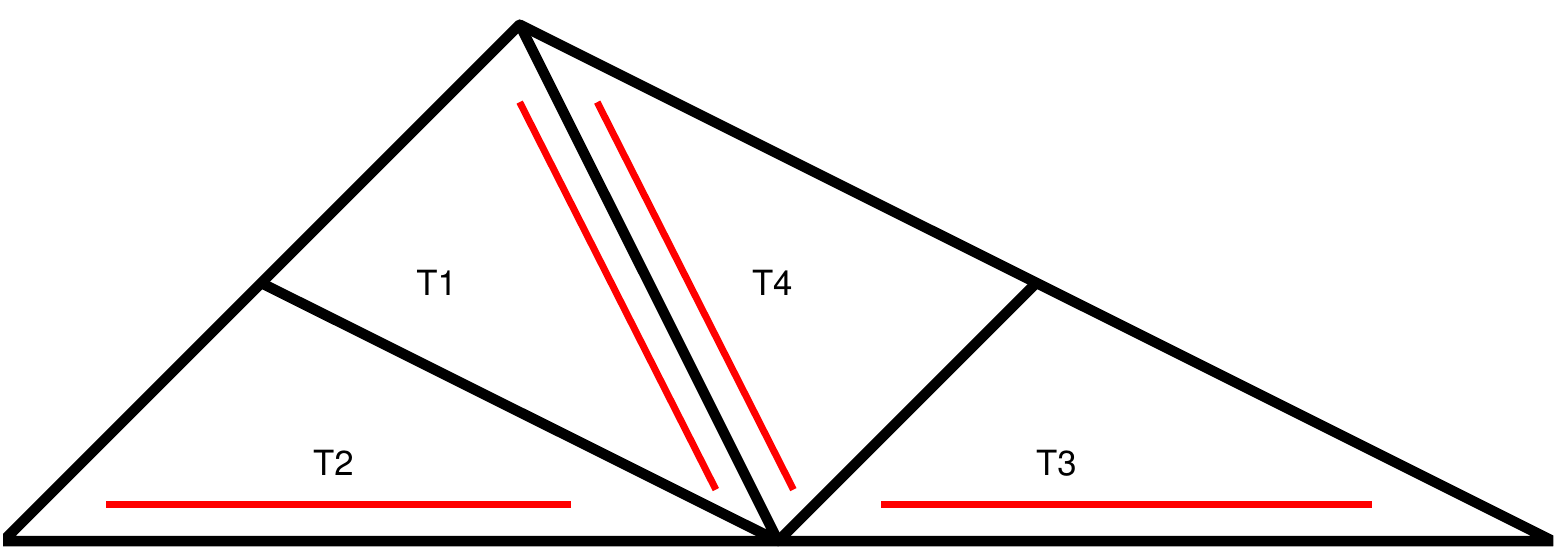}
 \caption{For~2D newest vertex bisection, each triangle $T\in\TT$ has one \emph{reference edge},
 indicated by the double line.
 Refinement of $T$ is done by bisecting
 the reference edge, where its midpoint becomes a new node.
 The reference edges of the son triangles are opposite to this newest vertex.
 Hanging nodes are avoided as follows:
 We  assume that certain edges of $T$, but at least the reference edge, are marked for refinement.
 By iterative bisection, refinement splits the element into 2, 3, or 4 son triangles, respectively.}
 \label{fig:nvb2d}
\end{figure}
\subsection{Mesh-refinement in~3D}
\label{section:weaksing:mesh3d}
For $d = 3$, a mesh $\TT_\coarse$ is a conforming triangulation of $\Gamma$ into compact plane surface triangles.
For mesh-refinement, we employ~2D newest vertex bisection (NVB);
see Figure~\ref{fig:nvb2d} and, e.g., \cite{stevenson,kpp}.
In particular, we assume that marked elements are bisected by three bisections into four son elements.
We assume that all meshes $\TT_\coarse \in \T = \refine(\TT_0)$
are obtained by applying this mesh-refinement strategy to a given initial mesh $\TT_0$.
As already mentioned in Section~\ref{section:twolevel}, this guarantees~\eqref{assumption:sons}--\eqref{assumption:closure}.
Moreover, NVB ensures that only finitely many shapes of triangles are generated.
In particular, all meshes $T_\coarse \in \T$ are uniformly $\gamma$-shape regular in the sense of
\begin{equation}\label{eq:2d:shape-regular}
 \max_{\TT_\coarse \in \T} \max_{T \in \TT_\coarse} \frac{{\rm diam}(T)^2}{{\rm area}(T)} \le \gamma < \infty,
\end{equation}
where $\gamma > 0$ depends only on $\TT_0$.
We note that conformity and~\eqref{eq:2d:shape-regular} also imply~\eqref{eq:1d:shape-regular}
(with a different constant $\gamma$ though).
We recall that $\widehat\TT_\coarse = \refine(\TT_\coarse; \TT_\coarse)$ denotes the uniform refinement,
where all elements have been bisected by three bisections. 
\end{subequations}

\subsection{Galerkin discretization}
\label{section:weaksing:galerkin}

Let $\TT_\coarse$ be a mesh. 
For the Galerkin discretization~\eqref{eq:abstract:galerkin} of~\eqref{eq:weaksing} resp.\ \eqref{eq:weaksing:helmholtz}, we consider the space of piecewise constant functions
\begin{equation}\label{eq:weaksing:Xh}
 \XX_\coarse :=
 \set{v_\coarse \in L^2(\Gamma)}{v_\coarse \vert_T \text{ is constant for all } T \in \TT_\coarse} \subset \H^{-1/2}(\Gamma).
\end{equation}
For $T \in \TT_\coarse$, let $\chi_T$ be the characteristic function, i.e., $\chi_T(x) = 1$ for $x \in T$ and $\chi_T(x) = 0$ for $x \in \Gamma\backslash T$. Then, $\set{\chi_T}{T \in \TT_\coarse}$ is the canonical basis of $\XX_\coarse$.

\begin{proposition}
The functional analytic framework of Section~\ref{section:weaksing:setting} and Section~\ref{section:weaksing:setting:helmholtz} together with its discretization in Section~\ref{section:weaksing:mesh2d}--\ref{section:weaksing:galerkin} satisfies all assumptions from Section~\ref{section:twolevel:abstract}. In particular, $\edual\cdot\cdot$ is always the scalar product~\eqref{eq:weaksing} induced by the Laplace single-layer operator~\eqref{eq:Vu=f}, while $\KK = 0$ for Laplace BEM (Section~\ref{section:weaksing:setting}) resp.\ $\KK = V_\kappa - V$ for Helmholtz BEM (Section~\ref{section:weaksing:setting:helmholtz}).
\end{proposition}%

\begin{proof}
It is well-known (see, e.g.,~\cite[Lemma~3.9.8]{sauter-schwab} or~\cite[Section~6.9]{steinbach}) that $\KK = V_\kappa - V : \H^{-1/2}(\Gamma) \to H^{1/2}(\Gamma)$ is a compact operator.
Hence,~\eqref{eq:weaksing} and~\eqref{eq:weaksing:helmholtz} meet the abstract variational formulation~\eqref{eq:abstract:weakform} of Section~\ref{section:twolevel:abstract}.
According to Remark~\ref{remark:approximation}, it only remains to verify that uniform mesh-refinement leads
to convergence of the best approximation error, i.e.,
that the definition $\TT^{(0)}_0 := \TT_0$ and $\TT^{(j+1)}_0 := \refine(\TT^{(j)}_0, \TT^{(j)}_0)$ for all $j \in \N_0$ guarantees that
\begin{equation*}
\lim_{j \to \infty} \min_{v_0^{(j)} \in \XX_0^{(j)}} \enorm{v - v_0^{(j)}} = 0
 \quad \text{for all } v \in \HH.
\end{equation*}
For $v \in L^2(\Gamma)$, we recall that 
\begin{equation*}
  \min_{v_0^{(j)} \in \XX_0^{(j)}} \enorm{v - v_0^{(j)}}
  \lesssim \norm{h_0^{(j)} v}{L^2(\Gamma)} \to 0
  \quad \text{as } j \to \infty.
\end{equation*}
Since $L^2(\Gamma)$ is dense in $\HH = \H^{-1/2}(\Gamma)$, this concludes the proof.
\end{proof}

\begin{figure}[t]
\def\ph{\widehat\varphi}
\def\plot#1#2{%
 \begin{minipage}{50mm}\begin{center}\includegraphics[width=40mm,height=20mm]{figures/#1.pdf}\\\scriptsize#2\end{center}\end{minipage}%
}
\plot{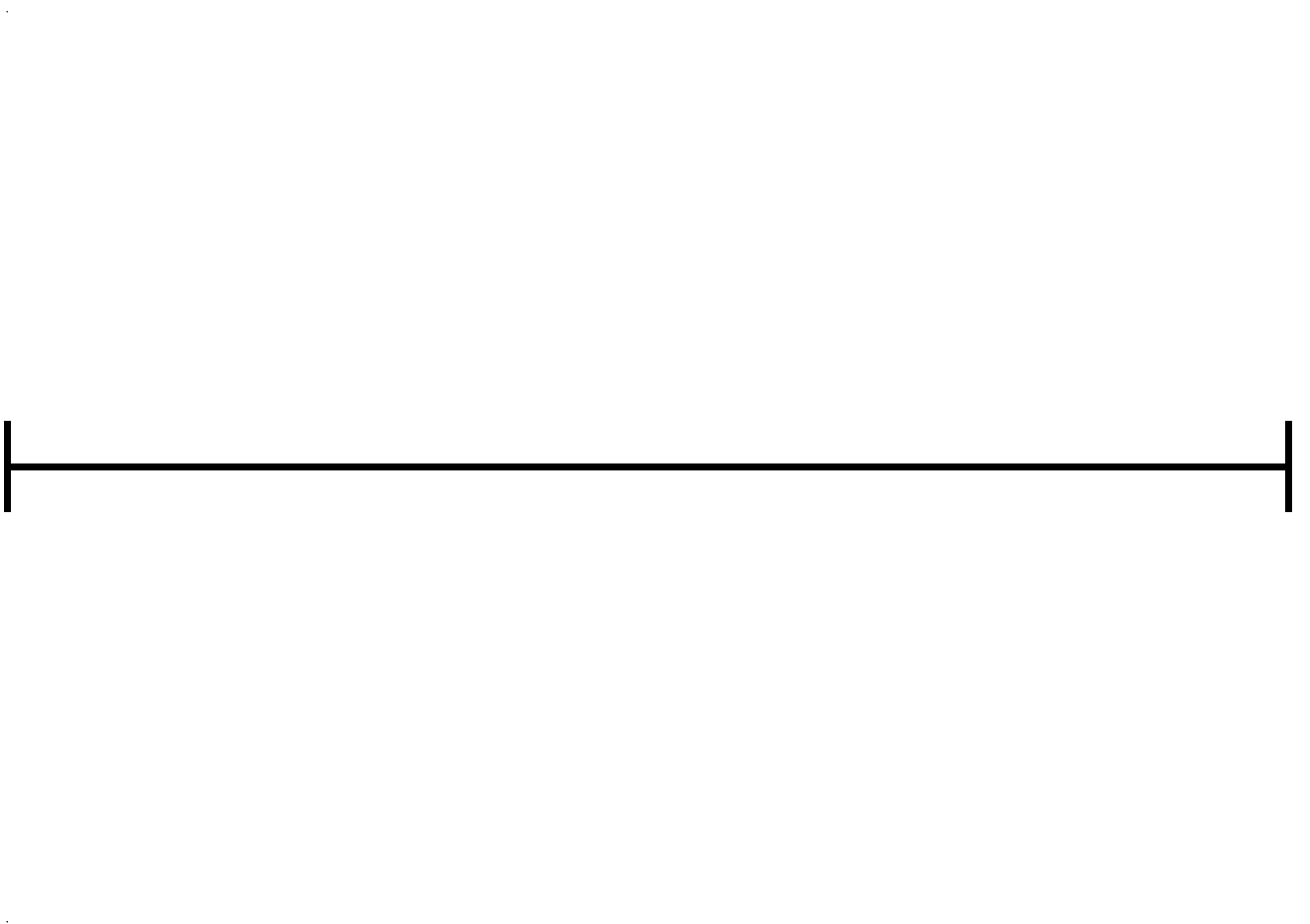}{Element $T \in \TT_\coarse$}
\plot{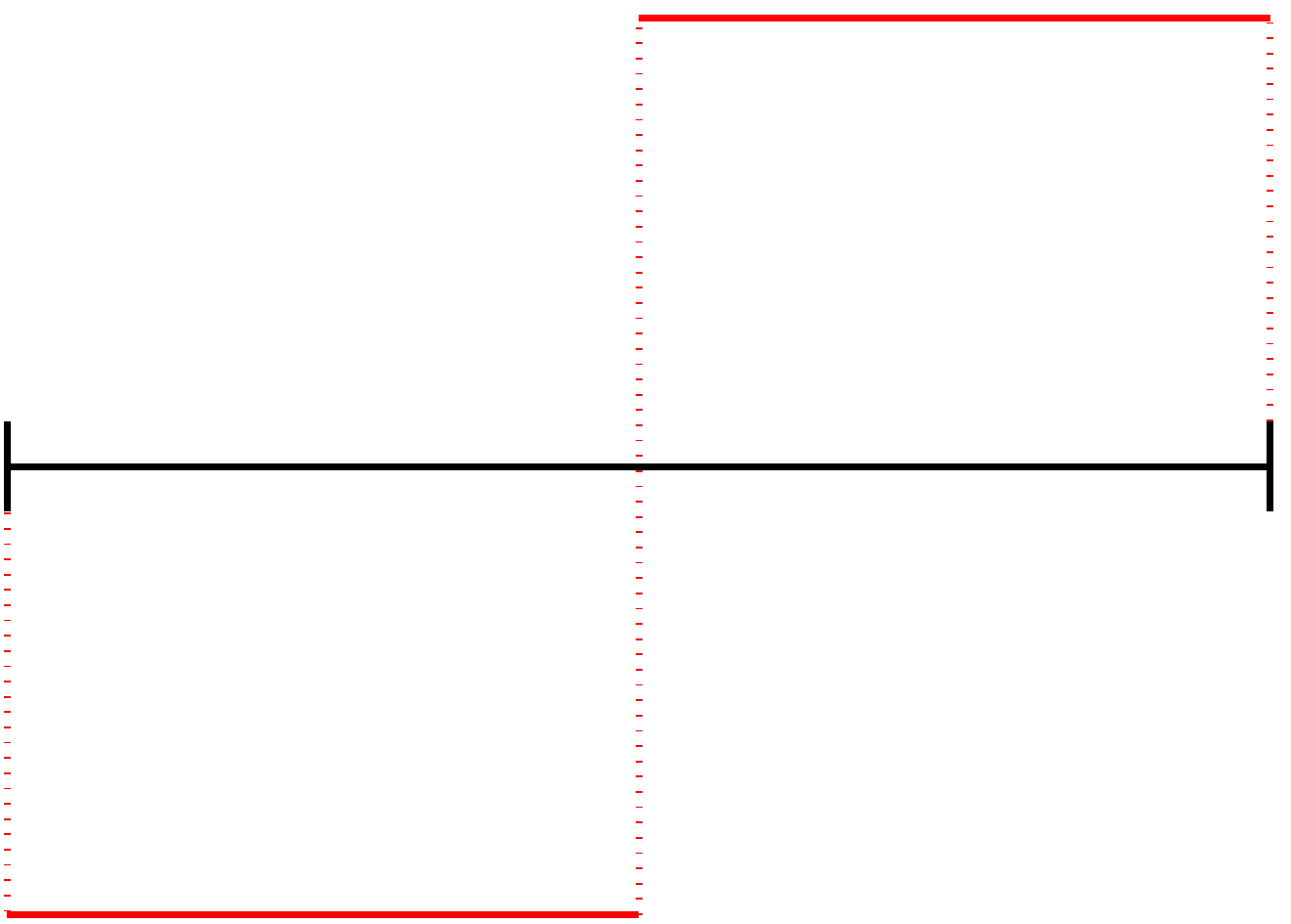}{$\ph_{\coarse,T}^{\,1} \in \widehat\XX_\coarse$.}
\caption{To generate the uniform refinement $\widehat\TT_\coarse$ for $d = 2$,
each coarse-mesh element $T \in \TT_\coarse$ (left) is bisected into two son elements of half length.
For the two-level estimator, $T$ is associated with one fine-mesh function $\ph_{\coarse,T}^{\,1} \in \widehat\XX_\coarse$
with $\supp(\ph_{\coarse,T}^{\,1}) = T$,
which takes the values $\pm1$ and is $L^2(T)$-orthogonal to the characteristic function $\chi_T$ (right).}
\label{fig:varphi2d}
\end{figure}

\begin{figure}[h]
\def\ph{\widehat\varphi}
\def\plot#1#2#3#4#5#6{%
 \psfrag{T0}{}%
 \psfrag{T1}{\tiny#1}%
 \psfrag{T2}{\tiny#2}%
 \psfrag{T3}{\tiny#3}%
 \psfrag{T4}{\tiny#4}%
 \begin{minipage}{30mm}\centering
 \includegraphics[width=\textwidth]{figures/#5.pdf}
 \\\scriptsize#6%
 \end{minipage}%
}
\begin{center}%
\plot{$+$}{$+$}{$+$}{$+$}{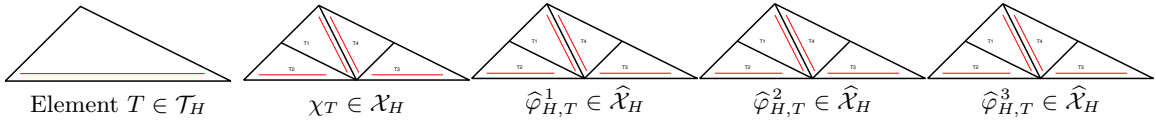}{Element $T \in \TT_\coarse$}
\plot{$+$}{$+$}{$+$}{$+$}{bisec3_01}{$\chi_T \in \XX_\coarse$}
\plot{$+$}{$+$}{$-$}{$-$}{bisec3_01}{$\ph_{\coarse,T}^{\,1} \in \widehat\XX_\coarse$}
\plot{$+$}{$-$}{$-$}{$+$}{bisec3_01}{$\ph_{\coarse,T}^{\,2} \in \widehat\XX_\coarse$}
\plot{$+$}{$-$}{$+$}{$-$}{bisec3_01}{$\ph_{\coarse,T}^{\,3} \in \widehat\XX_\coarse$}
\end{center}
 \caption{To get the uniform refinement $\widehat\TT_\coarse$ for $d = 3$, each coarse-mesh element $T \in \TT_\coarse$ is refined into four sons using three bisections. For the two-level estimator, $T$ is associated with $D = 3$ fine-mesh functions $\ph_{\coarse,T}^{\,j} \in \widehat\XX_\coarse$
 with $\supp(\ph_{\coarse,T}^{\,j}) = T$, which take the values $\pm1$ and are pairwise $L^2(T)$-orthogonal (also to the characteristic function $\chi_T$).}
 \label{fig:varphi3d}
\end{figure}

\subsection{Two-level \textsl{a~posteriori} error estimation}

For the weakly-singular integral equation~\eqref{eq:Vu=f} and lowest-order BEM~\eqref{eq:weaksing:Xh},
the local contributions of the two-level error estimator~\eqref{eq:twolevel} read
\begin{equation}\label{eq1:weaksing:twolevel}
 \tau_{\coarse}(T)^2 
 := \sum_{j = 1}^D \tau_{\coarse}(\widehat\varphi_{\coarse,T}^{\,j})^2
 \quad \text{with} \quad 
 \tau_{\coarse}(\widehat\varphi_{\coarse,T}^{\,j}) :=\frac{|\dual{f}{\widehat\varphi_{\coarse,T}^{\,j}} - b(u_\coarse,\widehat\varphi_{\coarse,T}^{\,j})|}{\enorm{\widehat\varphi_{\coarse,T}^{\,j}}},
\end{equation}
where $\{\widehat\varphi_{\coarse,T}^{\,j} \}_{j=1}^D$ is the set of local fine-mesh functions
with $D = 1$ for $d = 2$ (see Figure~\ref{fig:varphi2d}) and $D=3$ for $d = 3$ (see Figure~\ref{fig:varphi3d}).
Firstly considered in~\cite{mms1997,msw1998} for error control on quasi-uniform meshes,
we refer, e.g., to~\cite{eh2006,effp2009} for $h$-refinement in~2D and~3D
and to~\cite{hms2001} for $hp$-refinement in~2D for the fact that
\begin{equation}\label{eq2:weaksing:twolevel}
 \Ctwo^{-1} \, \enorm{\widehat u_\coarse - u_\coarse} \le \tau_\coarse \le \Ctwo \, \enorm{\widehat u_\coarse - u_\coarse},
\end{equation} 
where $\Ctwo \ge 1$ depends only on $\gamma$-shape regularity of $\widehat\TT_\coarse$ and hence only on $\TT_0$. The main observation for the proof of~\eqref{eq2:weaksing:twolevel} is that the decomposition
\begin{equation}
 \widehat\XX_\coarse = \XX_\coarse
 \oplus {\rm span}\set{\widehat\varphi_{\coarse,T}^{\, j}}{T \in \TT_\coarse, \, j = 1,\dots, D}
\end{equation}
is stable with respect to the $\H^{-1/2}$-norm.

\subsection{Convergence of Algorithm~\ref{algorithm} for weakly-singular integral equation}

The following theorem states that the weakly-singular integral equation is covered by the abstract framework of Section~\ref{section:twolevel}.

\begin{theorem}[optimal adaptivity for weakly-singular integral equation]\label{proposition:weaksing:twolevel}
Let $\TT_\coarse$ be a given mesh, $\MM_\coarse \subseteq \TT_\coarse$, and $\TT_\fine := \refine(\TT_\coarse,\MM_\coarse)$.
Then, the two-level estimator $\tau_\coarse^2 := \sum_{T \in \TT_\coarse}  \tau_{\coarse}(T)^2$ with its local contributions~\eqref{eq1:weaksing:twolevel} satisfies~\eqref{assumption:reliable}--\eqref{assumption:stable}, 
where $\Ctwo, \Cstab > 0$  depend only on $\gamma$-shape regularity of $\TT_\fine$ and hence on $\TT_0$. 
Applied to the framework of the weakly-singular integral equation for the Laplacian (Section~\ref{section:weaksing:setting}) or the Helmholtz problem (Section~\ref{section:weaksing:setting:helmholtz}), Algorithm~\ref{algorithm} thus leads to linear convergence~\eqref{eq1:linear} of the energy error with optimal algebraic rates~\eqref{eq:optimal}, if the saturation assumption~\eqref{eq:saturation} is satisfied.
\end{theorem}

\begin{proof}
Note that the proofs in~\cite{eh2006,effp2009} imply also that the decomposition
\begin{equation}\label{eq:weaksing:decomposition}
 \XX_\fine = \XX_\coarse
 \oplus {\rm span}\set{\widehat\varphi_{\coarse,T}^{\, j} \in \XX_\fine}{T \in \TT_\coarse, \, j = 1,\dots, D}
\end{equation}
is stable with respect to the $\H^{-1/2}$-norm. Arguing as in~\cite{eh2006,effp2009}, this leads to
\begin{equation*}
 \enorm{u_\fine - u_\coarse}^2
 \simeq \sum_{T \in \TT_\coarse} \sum_{\substack{j = 1 \\ \widehat\varphi_{\coarse,T}^{\,j} \in \XX_\fine}}^D \tau_{\coarse}(\widehat\varphi_{\coarse,T}^{\,j})^2,
\end{equation*}
where the two-level estimator on the right-hand side now involves only the enrichment from $\XX_\coarse$ to $\XX_\fine$; see~\eqref{eq:weaksing:decomposition}.
For details, we refer, e.g., to~\cite[Lemma~4.4, Proposition~4.5, Theorem~4.6]{effp2009} and to~\cite{mms1997} for the Helmholtz problem.
Next, note that
\begin{align*}
 \set{\widehat\varphi_{\coarse,T}^{\, j}}{T \in \MM_\coarse, \, j = 1,\dots, D}
 &\subseteq 
 \set{\widehat\varphi_{\coarse,T}^{\, j} \in \XX_\fine}{T \in \TT_\coarse, \, j = 1,\dots, D}
 \\
 &\subseteq 
 \set{\widehat\varphi_{\coarse,T}^{\, j}}{T \in \TT_\coarse \backslash\TT_\fine, \, j = 1,\dots, D},
\end{align*}
since the sons of marked elements $T \in \TT_\coarse$ are the same in $\TT_\fine$ and $\widehat\TT_\coarse$.
Clearly, this implies~\eqref{assumption:reliable}.

To prove~\eqref{assumption:stable}, note that the functions $\widehat\varphi_{\coarse,T}^{\, j}$ depend only on $T$, but not on the mesh $\TT_\coarse$.
Hence, the triangle inequality proves that
\begin{align*}
 &\big| \tau_\fine(\TT_\coarse \cap \TT_\fine) - \tau_\coarse(\TT_\coarse \cap \TT_\fine) \big|
 = \bigg| \Big( \sum_{T \in\TT_\coarse \cap \TT_\fine} \tau_\fine(T)^2 \Big)^{1/2} 
 - \Big( \sum_{T \in\TT_\coarse \cap \TT_\fine} \tau_\coarse(T)^2 \Big)^{1/2} \bigg|
 \\& \qquad
 \le \Big( \sum_{T \in\TT_\coarse \cap \TT_\fine} \sum_{j = 1}^D  \frac{|b(u_\fine - u_\coarse, \widehat\varphi_{\coarse,T}^{\,j})|^2}{\enorm{\widehat\varphi_{\coarse,T}^{\,j}}^2} \Big)^{1/2}
 \lesssim \enorm{u_\fine - u_\coarse},
\end{align*}
where the final estimate follows as in~\cite{eh2006,effp2009}. This concludes the proof.
\end{proof}

Even without the saturation assumption, we can adapt some ideas from~\cite{msv2008} to guarantee that Algorithm~\ref{algorithm} ensures at least convergence of the two-level estimator.

\begin{theorem}[plain convergence without the saturation assumption]\label{theorem:weaksing:plain}
For any adaptivity parameters $0 < \theta \le 1$ and $1 \le \Cmark \le \infty$ and independently of any saturation assumption, Algorithm~\ref{algorithm} guarantees that the two-level estimator $\tau_\ell^2 := \sum_{T \in \TT_\ell}  \tau_\ell(T)^2$ with its local contributions~\eqref{eq1:weaksing:twolevel} satisfies  $\tau_\ell \to 0$ as $\ell \to \infty$.
\end{theorem}

\def\est{\tau}
\begin{proof}
The proof is split into five steps.

{\bf Step~1.} There are several definitions of the fractional-order Sobolev space $H^{1/2}(\omega)$ for $\omega \subseteq \Gamma$,
e.g., by lifting to $H^1(\Omega)$, by use of the Sobolev--Slobodeckij seminorm, or by (real or complex) interpolation.
While these definitions lead to the same space, the norms are only equivalent up to constants,
which depend on $\omega$.
For $\omega \neq \Gamma$, the definition thus matters.
In what follows, we define $H^{1/2}(\omega) := [L^2(\omega); H^1(\omega)]_{1/2}$ by the K-method of interpolation~\cite[Section~3.1]{bergh}. Moreover, $\H^{-1/2}(\omega) := H^{1/2}(\omega)'$ denotes the dual space with respect to the extended $L^2(\omega)$ scalar product.

{\bf Step~2.} Let $\omega_\ell \subseteq \Gamma$ be a sequence of measurable subsets with $|\omega_\ell| \to 0$ as $\ell \to \infty$. The no-concentration of Lebesgue functions then implies that $\norm{v}{L^2(\omega_\ell)} \to 0$ as $\ell \to \infty$ for all $v \in L^2(\Gamma)$. Consequently, for $v \in H^1(\Gamma)$, it also follows that $\norm{v}{H^1(\omega_\ell)} \to 0$ as $\ell \to \infty$, and the interpolation estimate reveals that
\begin{equation*}
 \norm{v}{H^{1/2}(\omega_\ell)} 
 \lesssim \norm{v}{L^2(\omega_\ell)}^{1/2} \norm{v}{H^1(\omega_\ell)}^{1/2} 
 \xrightarrow{\ell \to \infty}0 \quad \text{for all } v \in H^1(\Gamma). 
\end{equation*}
Since $H^1(\Gamma)$ is dense in $H^{1/2}(\Gamma)$, it follows that 
\begin{equation*}
 \norm{v}{H^{1/2}(\omega_\ell)} \xrightarrow{\ell \to \infty}0 
 \quad \text{for all } v \in H^{1/2}(\Gamma).
\end{equation*}%

{\bf Step~3.} Define $H^0(T) := L^2(T)$. Let $\UU_\ell \subseteq \TT_\ell$. Since the estimate
\begin{equation}\label{eq:weaksing:step3}
 \sum_{T \in \UU_\ell} \norm{v}{H^s(T)}^2
 \le \norm{v}{H^s(\bigcup_{T \in \UU_\ell} T)}^2
 \quad \text{for all } v \in H^s(\Gamma)
\end{equation}
holds for $s \in \{0,1\}$, interpolation theory~\cite{bergh} implies that it also holds for $s = 1/2$.

{\bf Step~4.} Recall that the enrichtments $\widehat\varphi_{\ell,T}^{\,j}$ satisfy that $\dual{\widehat\varphi_{\ell,T}^{\,j}}{\chi_T} = 0$. Therefore, it follows from interpolation of the Poincar\'e inequality and a duality argument that
\begin{equation*}
 \norm{\widehat\varphi_{\ell,T}^{\,j}}{\H^{-1/2}(T)}
 \lesssim \diam(T)^{1/2} \, \norm{\widehat\varphi_{\ell,T}^{\,j}}{L^2(T)}
 \le \norm{h_\ell^{1/2} \, \widehat\varphi_{\ell,T}^{\,j}}{L^2(T)},
\end{equation*}
where $h_\ell \in \PP^0(\TT_\ell)$ denotes the local mesh-width function defined by $h_\ell|_{T'} := \diam(T')$; see, e.g.,~\cite[Theorem~4.1]{cp2006}. Together with an inverse estimate from~\cite[Theorem~3.6]{ghs2005}, this leads to
\begin{equation*}
 \norm{\widehat\varphi_{\ell,T}^{\,j}}{\H^{-1/2}(T)}
 \lesssim \norm{h_\ell^{1/2} \, \widehat\varphi_{\ell,T}^{\,j}}{L^2(T)}
 \lesssim \norm{\widehat\varphi_{\ell,T}^{\,j}}{\H^{-1/2}(\Gamma)}
 \simeq \enorm{\widehat\varphi_{\ell,T}^{\,j}}.
\end{equation*}%

{\bf Step~5.} With the aforegoing observations, the proof of the theorem essentially follows the lines of the proof of~\cite[Theorem~2.1]{msv2008}:
\begin{itemize}
\item the variational formulation~\eqref{eq:weaksing} clearly fits into the class of problems considered in~\cite[Section~2.1]{msv2008};
\item the Galerkin discretization~\eqref{eq:abstract:galerkin} with the discrete spaces~\eqref{eq:weaksing:Xh} satisfies the assumptions in~\cite[eq.~(2.6)--(2.8)]{msv2008};
\item the refinement strategies from Section~\ref{section:weaksing:mesh2d}--\ref{section:weaksing:mesh3d} satisfy the assumptions in~\cite[eq.~(25) and~(2.14)]{msv2008};
\end{itemize}
Let $\XX_\infty := \overline{\bigcup_{\ell=0}^\infty} \XX_\ell$ be the ``discrete limit space'', where the closure is understood with respect to $\HH = \H^{-1/2}(\Gamma)$. According to~\cite[Lemma~4.2]{msv2008}, there exists a unique $u_\infty \in \XX_\infty$ such that
\begin{equation}\label{eq:u_infty}
 b(u_\infty, v_\infty) = \dual{f}{v_\infty} 
 \quad \text{for all } v_\infty \in \XX_\infty,
\end{equation}
and it holds that $\enorm{u_\infty - u_\ell} \to 0$ as $\ell \to \infty$. Let $\TT_\infty := \bigcup_{k \ge 0} \bigcap_{\ell \ge k} \TT_\ell$ be the set of all elements which remain unrefined after finitely many steps of refinement. In the spirit of~\cite[eqs.~(4.10)]{msv2008}, for all $\ell \in \N_0$, we consider the decomposition
\begin{equation}\label{eq:T:decomposition}
 \TT_\ell = \TT_\ell^{\rm good} \cup \TT_\ell^{\rm bad} \cup  \TT_\ell^{\rm neither},
\end{equation}
where
\begin{align*}
 \TT_\ell^{\rm good} &:= \{T \in \TT_\ell : \widehat\varphi_{\ell,T}^j \in \XX_\infty \text{ for all } j = 1,\dots, D\} , \\
 \TT_\ell^{\rm bad} &:= \{T \in \TT_\ell : T' \in \TT_\infty \text{ for all } T' \in \TT_\ell \text{ with } T \cap T' \neq \emptyset \}, \\
 \TT_\ell^{\rm neither} &:= \TT_\ell \setminus (\TT_\ell^{\rm good} \cup \TT_\ell^{\rm bad}).
\end{align*}
In the following, we elaborate the ideas for the Laplacian (Section~\ref{section:weaksing:setting}). Replacing $V$ by $V_\kappa$, the same arguments apply to the Helmholtz problem (Section~\ref{section:weaksing:setting:helmholtz}).

The elements in $\TT_\ell^{\rm good}$ are refined sufficiently many times to guarantee that
\begin{equation}\label{eq2b:msv}
 \begin{split}
 &\tau_\ell(T)^2 
 \reff{eq1:weaksing:twolevel}= \sum_{j=1}^D \frac{|\dual{f - V u_\ell}{\widehat\varphi_{\ell,T}^j}|^2}{\enorm{\widehat\varphi_{\ell,T}^j}^2}
 \reff{eq:u_infty}= \sum_{j=1}^D \frac{|\dual{V(u_\infty - u_\ell)}{\widehat\varphi_{\ell,T}^j}|^2}{\enorm{\widehat\varphi_{\ell,T}^j}^2}
 \\& \quad
 \le \sum_{j=1}^D \norm{V(u_\infty - u_\ell)}{H^{1/2}(T)}^{2} \, \frac{\norm{\widehat\varphi_{\ell,T}^j}{\H^{-1/2}(T)}^2}{\enorm{\widehat\varphi_{\ell,T}^j}^2}
 \lesssim \norm{V(u_\infty - u_\ell)}{H^{1/2}(T)}^{2}.
 \end{split}
\end{equation}
The set $\TT_\ell^{\rm bad}$ consists of all elements such that the whole element patch remains unrefined.
The remaining elements are collected in the set $\TT_\ell^{\rm neither}$.
We note that $\TT_\ell^{\rm good}$ corresponds to $\GG_\ell^0$ in~\cite[eq.~(4.10a)]{msv2008} (but actually is a bit larger),
while $\TT_\ell^{\rm bad}$ coincides with the corresponding set $\GG_\ell^+$ in~\cite[eq.~(4.10b)]{msv2008}.
As a consequence, $\TT_\ell^{\rm neither}$ corresponds to $\GG_\ell^*$ in~\cite[eq.~(4.10c)]{msv2008}, but actually is a bit smaller.

With the mapping properties of $V$, we obtain that
\begin{equation}\label{eq:msv:step1}
 \begin{split}
 &\sum_{T \in \TT_\ell^{\rm good}} \tau_\ell(T)^2 
 \reff{eq2b:msv}{\lesssim} \sum_{T \in \TT_\ell^{\rm good}} \norm{V(u_\infty - u_\ell)}{H^{1/2}(T)}^2
 \reff{eq:weaksing:step3}\lesssim \norm{V(u_\infty - u_\ell)}{H^{1/2}(\Gamma)}^2
 \\& \qquad
 \simeq \norm{u_\infty - u_\ell}{\H^{-1/2}(\Gamma)}^2
 \simeq \enorm{u_\infty - u_{\ell}}^2 
 \xrightarrow{\ell \to \infty} 0.
\end{split}
\end{equation}
Let $\Gamma_\ell^{\rm neither} := \bigcup_{T \in \TT_\ell^{\rm neither}} T$. Since $\TT_\ell^{\rm neither}$ is contained in the corresponding set $\GG_\ell^*$ in~\cite[eq.~(4.10c)]{msv2008}, we may argue as in Step~1 of the proof of~\cite[Proposition~4.2]{msv2008} to show that $|\Gamma_\ell^{\rm neither}| \to 0$ as $\ell \to \infty$. According to Step~2, this leads to
\begin{align}\label{eq:msv:step2}
 &\sum_{T \in \TT_\ell^{\rm neither}} \tau_\ell(T)^2 
 \reff{eq1:weaksing:twolevel}= \sum_{T \in \TT_\ell^{\rm neither}} \sum_{j=1}^D \frac{|\dual{f - V u_\ell}{\widehat\varphi_{\ell,T}^j}|^2}{\enorm{\widehat\varphi_{\ell,T}^j}^2}
 \lesssim \sum_{T \in \TT_\ell^{\rm neither}}  \norm{f - V u_\ell}{H^{1/2}(T)}  
 \nonumber\\& \qquad 
 \reff{eq:weaksing:step3}\lesssim \norm{f - V u_\ell}{H^{1/2}(\Gamma_\ell^{\rm neither})}^2
 \lesssim\norm{f - V u_\infty}{H^{1/2}(\Gamma_\ell^{\rm neither})}^2
 + \norm{V(u_\infty - u_\ell)}{H^{1/2}(\Gamma_\ell^{\rm neither})}^2
 \nonumber\\& \qquad 
 \reff{eq:weaksing:step3}\lesssim\norm{f - V u_\infty}{H^{1/2}(\Gamma_\ell^{\rm neither})}^2
 + \norm{V(u_\infty - u_\ell)}{H^{1/2}(\Gamma)}^2
 \nonumber\\& \qquad 
 \lesssim\norm{f - V u_\infty}{H^{1/2}(\Gamma_\ell^{\rm neither})}^2
 + \enorm{u_\infty - u_\ell}^2 
 \xrightarrow{\ell \to \infty} 0. 
\end{align}%
Since $\MM_\ell \subseteq \TT_\ell \backslash \TT_\ell^{\rm bad} = \TT_\ell^{\rm good} \cup \TT_\ell^{\rm neither}$,
it follows from~\eqref{eq:msv:step1}--\eqref{eq:msv:step2} that
\begin{equation*}
 \theta \,\tau_\ell(T)^2
 \reff{eq:doerfler}\le \sum_{T \in \MM_\ell} \tau_\ell(T)^2
 \le \sum_{T \in \TT_\ell^{\rm good}} \tau_\ell(T)^2
 + \sum_{T \in \TT_\ell^{\rm neither}} \tau_\ell(T)^2
 \xrightarrow{\ell \to \infty} 0.
\end{equation*}
This concludes the proof.
\end{proof}

\begin{remark}
The convergence proof of Theorem~\ref{theorem:weaksing:plain} covers also other marking criterions than D\"orfler marking~\eqref{eq:doerfler}; for instance, the maximum criterion or estimator equilibration. To this end, only the final step in the proof of Theorem~\ref{theorem:weaksing:plain} has to be modified, where one can adapt the ideas of~\cite{msv2008}.
\end{remark}

\section{Hypersingular integral equation for $\boldsymbol{\Gamma \subsetneqq \partial\Omega}$}
\label{section:hypsing}

\subsection{Functional analytic framework for Laplace BEM}
\label{section:hypsing:setting}

Suppose the assumptions of Section~\ref{section:weaksing:setting}. In addition, suppose that the (relative) boundary of $\Gamma$ is non-trivial, i.e., $\partial\Gamma \neq \emptyset$, and hence $\Gamma$ is \emph{not} closed.
We consider the hypersingular integral equation
\begin{equation}\label{eq:Wu=f}
 (W u) (x) := -\partial_{\normal(x)}\int_\Gamma u(y) \, \partial_{\normal(y)}G(x-y) \d{y} = f(x)
 \quad \text{for } x \in \Gamma,
\end{equation}
where $\normal(\cdot)$ denotes the outer unit normal vector of $\partial\Omega$, $f \in H^{-1/2}(\Gamma)$ is some given right-hand side, and $u \in \H^{1/2}(\Gamma)$ is the sought integral density.
Here, we employ the notation $\H^{1/2}(\Gamma) := \set{v \in H^{1/2}(\partial\Omega)}{\supp(v) \subseteq \overline\Gamma}$ and $H^{-1/2}(\Gamma) := \H^{1/2}(\Gamma)'$ is its dual space with respect to the (extended) $L^2(\Gamma)$-scalar product. 

The variational formulation~\eqref{eq:abstract:weakform} of~\eqref{eq:Wu=f} reads:
Find $u \in \HH := \H^{1/2}(\Gamma)$ such that
\begin{equation}\label{eq:hypsing}
 b(u,v) := \edual{u}{v} := \dual{Wu}{v} = \dual{f}{v}
 \quad \text{for all } v \in \HH.
\end{equation}
It is well-known that the hypersingular integral operator $W: \H^{1/2}(\Gamma) \to H^{-1/2}(\Gamma)$ is a symmetric and elliptic isomorphism.
Therefore, $\edual\cdot\cdot$ is a scalar product, and the induced norm $\enorm{v}^2 := \edual{v}{v}$ is an equivalent norm on $\H^{1/2}(\Gamma)$. In particular, the Lax--Milgram theorem proves the existence and uniqueness of the solution $u \in \H^{1/2}(\Gamma)$ of~\eqref{eq:hypsing}.

\subsection{Functional analytic framework for Helmholtz BEM}
\label{section:hypsing:setting:helmholtz}

Recall the integral kernel $G_\kappa(z)$ from Section~\ref{section:weaksing:setting:helmholtz}. Adopting the notation from Section~\ref{section:hypsing:setting}, we consider the hypersingular integral equation
\begin{equation}\label{eq:Wu=f:hypsing}
 (W_\kappa u) (x) := -\partial_{\normal(x)}\int_\Gamma u(y) \, \partial_{\normal(y)}G_\kappa(x-y) \d{y} = f(x)
 \quad \text{for } x \in \Gamma,
\end{equation}
where $f \in H^{-1/2}(\Gamma)$ is some given right-hand side and $u \in \H^{1/2}(\Gamma)$ is the sought integral density.
The variational formulation~\eqref{eq:abstract:weakform} of~\eqref{eq:Wu=f:hypsing} reads:
Find $u \in \HH := \H^{1/2}(\Gamma)$ such that
\begin{equation}\label{eq:hypsing:helmholtz}
 b(u,v) := \dual{W_\kappa u}{v} = \dual{f}{v}
 \quad \text{for all } v \in \HH.
\end{equation}
It is well-known that the hypersingular operator $W_\kappa: \H^{1/2}(\Gamma) \to H^{-1/2}(\Gamma)$ is an isomorphism
if and only if $\kappa^2$ is not an eigenvalue of the interior Neumann problem for the Laplace operator;
see, e.g.,~\cite[Proposition~2.5]{steinbach}.
We suppose throughout that this is the case, i.e.,
$W_\kappa$ is an isomorphism and~\eqref{eq:hypsing:helmholtz} thus admits a unique solution $u \in \H^{1/2}(\Gamma)$.

\subsection{Galerkin discretization}
\label{section:hypsing:galerkin}

Let $\TT_\coarse$ be a mesh in the sense of
Sections~\ref{section:weaksing:mesh2d}--\ref{section:weaksing:mesh3d}.
For the Galerkin discretization~\eqref{eq:abstract:galerkin} of~\eqref{eq:hypsing} resp.\ \eqref{eq:hypsing:helmholtz}, we consider the space of continuous piecewise linear finite elements
\begin{equation}\label{eq:hypsing:galerkin}
 \XX_\coarse := \set{v_\coarse \in C(\overline\Gamma)}{v_\coarse \vert_T \text{ is affine for all } T \in \TT_\coarse \text{ and } v_\coarse|_{\partial\Gamma} = 0} \subset \H^{1/2}(\Gamma).
\end{equation}
Let $\NN_\coarse$ be the set of vertices of $\TT_\coarse$. 
For $z \in \NN_\coarse$, let $\varphi_{\coarse,z}$ be the associated hat function, i.e., $\varphi_{\coarse,z}$ is piecewise affine, globally continuous,
and satisfies the Kronecker property $\varphi_{\coarse,z}(z') = \delta_{zz'}$ for all $z' \in \NN_\coarse$.
Then, $\set{\varphi_{\coarse,z}}{z \in \NN_\coarse \setminus \partial \Gamma}$ is the standard basis of $\XX_\coarse$.

\begin{proposition}
The functional analytic framework of Section~\ref{section:hypsing:setting} and Section~\ref{section:hypsing:setting:helmholtz} together with its discretization in Section~\ref{section:hypsing:galerkin} satisfies all assumptions from Section~\ref{section:twolevel:abstract}. In particular, $\edual\cdot\cdot$ is always the scalar product~\eqref{eq:hypsing} induced by the Laplace hypersingular integral operator~\eqref{eq:Wu=f}, while $\KK = 0$ for Laplace BEM (Section~\ref{section:hypsing:setting}) resp.\ $\KK = W_\kappa - W$ for Helmholtz BEM (Section~\ref{section:hypsing:setting:helmholtz}).
\end{proposition}%

\begin{proof}
It is well-known (see, e.g.,~\cite[Lemma~3.9.8]{sauter-schwab} or~\cite[Section~6.9]{steinbach}) that $\KK = W_\kappa - W : \H^{-1/2}(\Gamma) \to H^{1/2}(\Gamma)$ is a compact operator. Hence,~\eqref{eq:hypsing} and~\eqref{eq:hypsing:helmholtz} meet the abstract variational formulation~\eqref{eq:abstract:weakform} of Section~\ref{section:twolevel:abstract}. According to Remark~\ref{remark:approximation}, it only remains to verify that uniform mesh-refinement leads to convergence of the best approximation error, i.e., that the definition $\TT^{(0)}_0 := \TT_0$ and $\TT^{(j+1)}_0 := \refine(\TT^{(j)}_0, \TT^{(j)}_0)$ for all $j \in \N_0$ guarantees that
\begin{equation*}
\lim_{j \to \infty} \min_{v_0^{(j)} \in \XX_0^{(j)}} \enorm{v - v_0^{(j)}} = 0
 \quad \text{for all } v \in \HH.
\end{equation*}
For $v \in \H^1(\Gamma)$, we recall that 
\begin{equation*}
  \min_{v_0^{(j)} \in \XX_0^{(j)}} \enorm{v - v_0^{(j)}}
  \lesssim \norm{h_0^{(j)} \nabla v}{L^2(\Gamma)} \to 0
  \quad \text{as } j \to \infty.
\end{equation*}
Since $\H^1(\Gamma)$ is dense in $\HH = \H^{1/2}(\Gamma)$, this concludes the proof.
\end{proof}

\subsection{Two-level \textsl{a~posteriori} error estimation}
Let $\TT_\coarse$ be a mesh with uniform refinement $\widehat\TT_\coarse$. 
In 2D, define $D = 1$ and let $\widehat\varphi_{\coarse,T}^1 \in \widehat\XX_\coarse$ be the hat function
associated with the midpoint of $T \in \TT_\coarse$.
In 3D, define $D = 3$.
For an element $T \in \TT_\coarse$ with its three edges $E_1, E_2, E_3$, let $\widehat\varphi_{\coarse,T}^{\,j} \in \widehat\XX_\coarse$ be either the zero function, if $E_j$ lies on the boundary $\partial\Gamma$ of $\Gamma$, or the hat function associated with the midpoint of $E_j$.

For the hypersingular integral equation~\eqref{eq:Wu=f} and lowest-order BEM~\eqref{eq:hypsing:galerkin}, the local contributions of the two-level error estimator~\eqref{eq:twolevel} read
\begin{equation}\label{eq1:hypsing:twolevel}
 \tau_{\coarse}(T)^2 
 := \sum_{j = 1}^D \tau_{\coarse}(\widehat\varphi_{\coarse,T}^{\,j})^2
 \quad \text{with} \quad 
 \tau_{\coarse}(\widehat\varphi_{\coarse,T}^{\,j}) :=\frac{|\dual{f}{\widehat\varphi_{\coarse,T}^{\,j}} - b(u_\coarse,\widehat\varphi_{\coarse,T}^{\,j})|}{\enorm{\widehat\varphi_{\coarse,T}^{\,j}}},
\end{equation}
where $D = 1$ for $d = 2$ and $D=3$ for $d = 3$.
Firstly considered in~\cite{mms1997,ms2000} for error control on quasi-uniform meshes, we refer, e.g., to~\cite{efgp2013,affkp2015} for $h$-refinement in~2D and~3D and to~\cite{hms2001,hs2001,heuer2002} for $hp$-refinement for the fact that
\begin{equation}\label{eq2:hypsing:twolevel}
 \Ctwo^{-1} \, \enorm{\widehat u_\coarse - u_\coarse} \le \tau_\coarse \le \Ctwo \, \enorm{\widehat u_\coarse - u_\coarse},
\end{equation} 
where $\Ctwo \ge 1$ depends only on $\gamma$-shape regularity of $\widehat\TT_\coarse$ and hence only on $\TT_0$.
The main observation for the proof of~\eqref{eq2:hypsing:twolevel} is that the decomposition
\begin{equation}\label{eq0:hypsing:decomposition}
 \widehat\XX_\coarse = \XX_\coarse
 \oplus {\rm span}\set{\widehat\varphi_{\coarse,z}}{z \in \widehat\NN_\coarse \backslash \NN_\coarse \text{ with } z \not\in\partial\Gamma}
\end{equation}
is stable with respect to the $\H^{1/2}$-norm. 

\subsection{Convergence of Algorithm~\ref{algorithm} for hypersingular integral equation}

The following theorem states that the hypersingular integral equation is also covered by the abstract framework of Section~\ref{section:twolevel}.

\begin{theorem}[optimal adaptivity for hypersingular integral equation]\label{proposition:hypsing:twolevel}
Let $\TT_\coarse$ be a given mesh, $\MM_\coarse \subseteq \TT_\coarse$, and $\TT_\fine := \refine(\TT_\coarse,\MM_\coarse)$.
Then, the two-level estimator $\tau_\coarse^2 := \sum_{T \in \TT_\coarse}  \tau_{\coarse}(T)^2$ with its local contributions~\eqref{eq1:hypsing:twolevel} satisfies~\eqref{assumption:reliable}--\eqref{assumption:stable}, 
where $\Ctwo, \Cstab > 0$  depend only on $\gamma$-shape regularity of $\TT_\fine$ and hence on $\TT_0$. 
Applied to the framework of the hypersingular integral equation for the Laplacian (Section~\ref{section:hypsing:setting}) or the Helmholtz problem (Section~\ref{section:hypsing:setting:helmholtz}), Algorithm~\ref{algorithm} thus leads to linear convergence~\eqref{eq1:linear} of the energy error with optimal algebraic rates~\eqref{eq:optimal}, if the saturation assumption~\eqref{eq:saturation} is satisfied.
\end{theorem}

\begin{proof}
The essential observation is that 
\begin{equation*}
 \widehat \varphi_{\coarse,T}^{\,j} \in \SS^1_0(\TT_h)
 \quad \text{if and only if the corresponding edge $E_j$ of $T \in \TT_\coarse$ is bisected.}
\end{equation*}
While this is clear for 2D, it is less obvious for 3D but a consequence of the refinement pattern of NVB; see Figure~\ref{fig:nvb2d}.
First, this allows to rewrite~\eqref{eq0:hypsing:decomposition} in the form
\begin{equation*}
 \widehat\XX_\coarse = \XX_\coarse
 \oplus {\rm span}\set{\widehat\varphi_{\coarse,T}^{\,j}}{T \in \TT_\coarse, \, j = 1,\dots, D}.
\end{equation*}
Second, with this observation, the proofs in~\cite{efgp2013,affkp2015} imply also that the decomposition
\begin{equation}\label{eq:hypsing:decomposition}
 \XX_\fine = \XX_\coarse
 \oplus {\rm span}\set{\widehat\varphi_{\coarse,T}^{\,j} \in \XX_\fine}{T \in \TT_\coarse, \, j = 1,\dots, D}
\end{equation}
is stable with respect to the $\H^{1/2}$-norm. Arguing as in~\cite{efgp2013,affkp2015}, this leads to
\begin{equation*}
 \enorm{u_\fine - u_\coarse}^2
 \simeq \sum_{T \in \TT_\coarse} \sum_{\substack{j = 1 \\ \widehat\varphi_{\coarse,T}^{\,j} \in \XX_\fine}}^D \tau_{\coarse}(\widehat\varphi_{\coarse,T}^{\,j})^2,
\end{equation*}
where the two-level estimator on the right-hand side now involves only the enrichment from $\XX_\coarse$ to $\XX_\fine$; see~\eqref{eq:hypsing:decomposition}.
For details, we refer, e.g., to~\cite[Lemma~4.4, Theorem~4.1]{efgp2013} and~\cite[Lemma~15, Lemma~16, Theorem~14]{affkp2015} and to~\cite{mms1997} for the Helmholtz problem.
Next, note that
\begin{align*}
 \set{\widehat\varphi_{\coarse,T}^{\,j}}{T \in \MM_\coarse, \, j = 1,\dots, D}
 &\subseteq 
 \set{\widehat\varphi_{\coarse,T}^{\,j} \in \XX_\fine}{T \in \TT_\coarse, \, j = 1,\dots, D}
 \\
 &\subseteq 
 \set{\widehat\varphi_{\coarse,T}^{\,j}}{T \in \TT_\coarse \backslash\TT_\fine, \, j = 1,\dots, D},
\end{align*}
since the sons of marked elements $T \in \TT_\coarse$ are the same in $\TT_\fine$ and $\widehat\TT_\coarse$.
Clearly, this implies~\eqref{assumption:reliable}.

To prove~\eqref{assumption:stable}, note that the functions $\widehat\varphi_{\coarse,T}^{\,j}$ depend only on $T$, but not on the mesh, i.e., for $T \in \TT_\coarse \cap \TT_\fine$, the hat function $\widehat\varphi_{\coarse,T}^{\,j}$ belongs to $\widehat\XX_\fine \cap \widehat\XX_\coarse$.
Hence, the triangle inequality proves that
\begin{align*}
 &\big| \tau_\fine(\TT_\coarse \cap \TT_\fine) - \tau_\coarse(\TT_\coarse \cap \TT_\fine) \big|
 = \bigg| \Big( \sum_{T \in\TT_\coarse \cap \TT_\fine} \tau_\fine(T)^2 \Big)^{1/2} 
 - \Big( \sum_{T \in\TT_\coarse \cap \TT_\fine} \tau_\coarse(T)^2 \Big)^{1/2} \bigg|
 \\& \qquad
 \le \Big( \sum_{T \in\TT_\coarse \cap \TT_\fine} \sum_{j = 1}^D  \frac{|\edual{u_\fine - u_\coarse}{\widehat\varphi_{\coarse,T}^{\,j}}|^2}{\enorm{\widehat\varphi_{\coarse,T}^{\,j}}^2} \Big)^{1/2}
 \lesssim \enorm{u_\fine - u_\coarse},
\end{align*}
where the final estimate follows as in~\cite{efgp2013,affkp2015}. This concludes the proof.
\end{proof}

Even without the saturation assumption, we can argue as in Section~\ref{section:weaksing}
(Theorem~\ref{theorem:weaksing:plain}) and show that Algorithm~\ref{algorithm} ensures at least
convergence of the two-level estimator.

\begin{theorem}[plain convergence without the saturation assumption]\label{theorem:hypsing:plain}
For any adaptivity parameters $0 < \theta \le 1$ and $1 \le \Cmark \le \infty$ and independently of any saturation assumption, Algorithm~\ref{algorithm} guarantees that the two-level estimator $\tau_\ell^2 := \sum_{T \in \TT_\ell}  \tau_{\ell}(T)^2$ with its local contributions~\eqref{eq1:hypsing:twolevel} satisfies  $\tau_\ell \to 0$ as $\ell \to \infty$.
\end{theorem}

\begin{proof}
The proof is split into five steps.
The first four steps prove some technical results on fractional-order Sobolev norms.
Even though they might be well-known to the experts, we include their brief proofs for the convenience of the reader;
see also the discussion in Step~1 of the proof of Theorem~\ref{theorem:weaksing:plain}
for the fact that the definition of the fractional norms matters.

{\bf Step~1.}
For $\omega \subseteq \Gamma$, we recall that $H^{-s}(\omega) := \H^s(\omega)'$,
where $\H^s(\omega) = \set{v \in H^s(\Gamma)}{\supp(v) \subseteq \overline{\omega}}$ for all $s \in [0,1]$.
For $s \in (0,1)$, the space is equipped with the interpolation norm.
For all $\omega, \widehat\omega \subseteq \Gamma$ with $\omega \subseteq \widehat\omega$,
we then have the norm equivalence
\begin{equation} \label{eq:hypsing:step1:equiv}
\norm{v}{\H^s(\omega)}
\simeq
\norm{v}{\H^s(\widehat\omega)}
\quad \text{for all }
v \in \H^{s}(\omega).
\end{equation}
Indeed, on the one hand, the estimate 
$\norm{v}{\H^s(\omega)} \leq \norm{v}{\H^s(\widehat\omega)}$
holds for $s \in \{0,1\}$,
and thus also for $s \in (0,1)$ by interpolation theory~\cite{bergh}.
On the other hand, the converse estimate follows, e.g., from the open mapping theorem,
observing that $\H^s(\omega)$ is a closed subspace of $\H^s(\widehat\omega)$.
As a consequence,
for all $\phi \in H^{-s}(\omega)$,
it follows that
\begin{equation} \label{eq:hypsing:step1:equiv2}
\norm{\phi}{H^{-s}(\omega)}
=
\sup_{v \in \H^s(\omega) \setminus \{0\}} \frac{\dual{\phi}{v}}{\norm{v}{\H^s(\omega)}}
\lesssim 
\sup_{v \in \H^s(\widehat\omega) \setminus \{0\}} \frac{\dual{\phi}{v}}{\norm{v}{\H^s(\widehat\omega)}}
=
\norm{\phi}{H^{-s}(\widehat\omega)}.
\end{equation}

{\bf Step~2.}
Let $\omega_\ell \subseteq \Gamma$ be a sequence of measurable subsets with $|\omega_\ell| \to 0$ as $\ell \to \infty$.
We show that $\norm{\phi}{H^{-1/2}(\omega_\ell)} \to 0$ as $\ell \to \infty$ for all $\phi \in H^{-1/2}(\Gamma)$.
To this end, let $\phi \in H^{-1/2}(\Gamma)$ and $\eps>0$ be arbitrary.
Since $L^2(\Gamma)$ is densely embedded in $H^{-1/2}(\Gamma)$,
there exists $\widetilde\phi \in L^2(\Gamma)$ such that $\norm{\phi - \widetilde\phi}{H^{-1/2}(\Gamma)} \leq \eps/2$.
By no-concentration of Lebesgue functions,
it holds that $\norm{\widetilde\phi}{L^2(\omega_\ell)} \to 0$ as $\ell \to \infty$.
In particular, there exists $\ell_0 \in \N_0$ such that
$\norm{\widetilde\phi}{L^2(\omega_\ell)} \leq \eps/2$ for all $\ell \geq \ell_0$. 
Hence, for all $\ell \geq \ell_0$, using also~\eqref{eq:hypsing:step1:equiv}, we obtain that
\begin{equation*}
\norm{\phi}{H^{-1/2}(\omega_\ell)}
\leq
\norm{\widetilde\phi}{H^{-1/2}(\omega_\ell)}
+
\norm{\phi - \widetilde\phi}{H^{-1/2}(\omega_\ell)}
\lesssim
\norm{\widetilde\phi}{L^2(\omega_\ell)}
+
\norm{\phi - \widetilde\phi}{H^{-1/2}(\Gamma)}
\leq \eps.
\end{equation*}
This proves the desired limit.

{\bf Step~3.}
We consider a partition $\omega = \bigcup_{m=1}^M \omega_m$ of a subset $\omega \subseteq \Gamma$
satisfying $\omega_n \cap \omega_m = \emptyset$ if $n \neq m$.
We show that, for all $s \in (0,1)$, it holds that
\begin{equation} \label{eq:auxiliary_coloring}
\sum_{m=1}^M \norm{\phi}{H^{-s}(\omega_m)}^2
\le
\norm{\phi}{H^{-s}(\omega)}^2
\quad
\text{for all }
\phi \in H^{-s}(\Gamma).
\end{equation}
To that end,
we consider the product space
$\Pi^s := \prod_{m=1}^M \H^s(\omega_m)$,
endowed with the product norm
$\norm{\boldsymbol{u}}{\Pi^s}^2 := \sum_{m=1}^M \norm{u_m}{\H^s(\omega_m)}$
for all $\boldsymbol{u} = (u_m)_{m=1,\dots,M}$.
We consider the sum operator (which coincides with the identity, because we are dealing with a partition of $\omega$)
$A^s: \Pi^s \to \H^{s}(\omega)$ defined by
$A^s\boldsymbol{u} := \sum_{m=1}^M u_m$ for all $\boldsymbol{u} = (u_m)_{m=1,\dots,M} \in \Pi^s$.
Since
\begin{equation*}
\norm{A^s\boldsymbol{u}}{\H^{s}(\omega)}^2
= \Big\lVert \sum_{m=1}^M u_m \Big\rVert_{\H^{s}(\omega)}^2
\le
\sum_{m=1}^M \norm{u_m}{\H^s(\omega_m)}
= \norm{\boldsymbol{u}}{\Pi^s}^2
\quad
\text{for all } \boldsymbol{u} = (u_m)_{m=1,\dots,M} \in \Pi^s
\end{equation*}
holds for $s \in \{ 0,1\}$ (even with equality sign), then the
inequality is true for all $s \in (0,1)$ by interpolation.
Hence, $A^s$ is linear and bounded for all $s \in (0,1)$ and the same holds for its adjoint
$(A^s)': \H^{s}(\omega)' \to (\Pi^s)'$.
Since $(\Pi^s)' = \big( \prod_{m=1}^M \H^s(\omega_m) \big)' = \prod_{m=1}^M \H^s(\omega_m)' = \prod_{m=1}^M H^{-s}(\omega_m)$, the estimate~\eqref{eq:auxiliary_coloring} follows from the boundedness of $(A^s)'$.

{\bf Step~4.}
For all $T \in \TT_\ell$, let $\omega_\ell(T) \subseteq \Gamma$ denote the element-patch of $T$ in $\TT_\ell$.
For all $\UU_\ell \subseteq \TT_\ell$ and $U_\ell := \bigcup_{T \in \UU_\ell} T$,
we define $\omega_\ell(U_\ell) := \bigcup_{T \in \UU_\ell} \omega_\ell(T)$.
We show that, for all $s \in (0,1)$, it holds that
\begin{equation} \label{eq:coloring}
\sum_{T \in \UU_\ell} \norm{\phi}{H^{-s}(\omega_\ell(T))}^2
\lesssim
\norm{\phi}{H^{-s}(\omega_\ell(U_\ell))}^2
\quad
\text{for all }
\phi \in H^{-s}(\Gamma),
\end{equation}
To this end, we note
that both the number of elements contained in each patch
and the number of patches to which each element belong
are uniformly bounded, with the bounds depending only on 
the shape-regularity of the mesh, and thus only on $\TT_0$.
With an inductive construction (see, e.g., \cite[Lemma~3.1]{cms2001}), we can construct a partition
$\UU_\ell = \bigcup_{m=1}^M \mathcal{Q}_m$ such that
$\mathcal{Q}_m \cap \mathcal{Q}_n = \emptyset$ if $m \neq n$
with the property that $\omega_\ell(T) \cap \omega_\ell(T') = \emptyset$
for all $T,T' \in \mathcal{Q}_m$ with $T \neq T'$ and $m=1,\dots,M$.
Again, $M$ depends only on the shape-regularity of $\TT_\ell$, and hence only on $\TT_0$.
Define $Q_m := \bigcup_{T \in \mathcal{Q}_m} T$.
For all $\phi \in H^{-s}(\Gamma)$, we then conclude that
\begin{equation*}
\sum_{T \in \UU_\ell} \norm{\phi}{H^{-s}(\omega_\ell(T))}^2
=
\sum_{m=1}^M
\sum_{T \in \mathcal{Q}_m}
\norm{\phi}{H^{-s}(\omega_\ell(T))}^2
\reff{eq:auxiliary_coloring}\le
\sum_{m=1}^M \norm{\phi}{H^{-s}(\omega_\ell(Q_m))}^2
\reff{eq:hypsing:step1:equiv2}\lesssim
M \norm{\phi}{H^{-s}(\omega_\ell(U_\ell))}^2.
\end{equation*}

{\bf Step~5.} 
We argue as for Step~5 of the the proof of Theorem~\ref{theorem:weaksing:plain}.
We sketch the proof for the Laplacian, but the same ideas apply to the Helmholtz problem
(Section~\ref{section:hypsing:setting:helmholtz}) replacing $W$ with $W_\kappa$.
Let $\XX_\infty := \overline{\bigcup_{\ell=0}^\infty} \XX_\ell$.
Recall from~\cite[Lemma~4.2]{msv2008} that there exists a unique $u_\infty \in \XX_\infty$ such that
\begin{equation}\label{eq:u_infty:hypsing}
 b(u_\infty, v_\infty) = \dual{f}{v_\infty} 
 \quad \text{for all } v_\infty \in \XX_\infty,
\end{equation}
and it holds that $\enorm{u_\infty - u_\ell} \to 0$ as $\ell \to \infty$.
Use the decomposition~\eqref{eq:T:decomposition}.
Arguing as for the weakly-singular integral equation, we see that
\begin{equation*}
\begin{split}
&\tau_\ell(T)^2 
\reff{eq1:hypsing:twolevel}= \sum_{j=1}^D \frac{|\dual{f - W u_\ell}{\widehat\varphi_{\ell,T}^j}|^2}{\enorm{\widehat\varphi_{\ell,T}^j}^2}
\reff{eq:u_infty:hypsing}= \sum_{j=1}^D \frac{|\dual{W(u_\infty - u_\ell)}{\widehat\varphi_{\ell,T}^j}|^2}{\enorm{\widehat\varphi_{\ell,T}^j}^2}
\\& \quad
\le \sum_{j=1}^D \norm{W(u_\infty - u_\ell)}{H^{-1/2}(\omega_\ell(T))} \, \frac{\norm{\widehat\varphi_{\ell,T}^j}{\H^{1/2}(\omega_\ell(T))}^2}{\enorm{\widehat\varphi_{\ell,T}^j}^2}
\lesssim \norm{W(u_\infty - u_\ell)}{H^{-1/2}(\omega_\ell(T))},
\end{split}
\end{equation*}
where we have used the norm equivalence
\begin{equation*}
\norm{\widehat\varphi_{\ell,T}^j}{\H^{1/2}(\omega_\ell(T))}
\reff{eq:hypsing:step1:equiv}\simeq
\norm{\widehat\varphi_{\ell,T}^j}{\H^{1/2}(\Gamma)}
\simeq
\enorm{\widehat\varphi_{\ell,T}^j}.
\end{equation*}
With Step~4 and the mapping properties of $W$, we obtain that
\begin{equation}\label{hypsing:eq:msv:step1}
 \begin{split}
 \sum_{T \in \TT_\ell^{\rm good}} \tau_\ell(T)^2 
 &
 \lesssim \sum_{T \in \TT_\ell^{\rm good}} \norm{W(u_\infty - u_\ell)}{H^{-1/2}(\omega_\ell(T))}^2
 \reff{eq:coloring}\lesssim \norm{W(u_\infty - u_\ell)}{H^{-1/2}(\omega_\ell( \Gamma_\ell^{\rm good} ))}^2
 \\&
\reff{eq:hypsing:step1:equiv2}\lesssim
\norm{W(u_\infty - u_\ell)}{H^{-1/2}(\Gamma)}^2
 \lesssim \norm{u_\infty - u_\ell}{\H^{1/2}(\Gamma)}^2
 \simeq \enorm{u_\infty - u_{\ell}}^2 
 \xrightarrow{\ell \to \infty} 0,
\end{split}
\end{equation}
where $\Gamma_\ell^{\rm good} := \bigcup_{T \in \TT_\ell^{\rm good}} T$.
Recall that $\Gamma_\ell^{\rm neither} := \bigcup_{T \in \TT_\ell^{\rm neither}} T$ satisfies that
$|\Gamma_\ell^{\rm neither}| \to 0$ as $\ell \to \infty$.
Moreover, it holds that 
$|\Gamma_\ell^{\rm neither}| \simeq |\omega_\ell(\Gamma_\ell^{\rm neither})|$.
By Step~2, this leads to
\begin{equation}\label{hypsing:eq:msv:step2}
 \begin{split}
& \sum_{T \in \TT_\ell^{\rm neither}} \tau_\ell(T)^2 
= \sum_{T \in \TT_\ell^{\rm neither}} \sum_{j=1}^D \frac{|\dual{f - W u_\ell}{\widehat\varphi_{\ell,T}^j}|^2}{\enorm{\widehat\varphi_{\ell,T}^j}^2}
 \\& \qquad
 \lesssim \sum_{T \in \TT_\ell^{\rm neither}}  \norm{f - W u_\ell}{H^{-1/2}(\omega_\ell(T))}^2
 \reff{eq:coloring}\lesssim \norm{f - W u_\ell}{H^{-1/2}( { \omega_\ell( \Gamma_\ell^{\rm neither} )})}^2
\to
0. 
 \end{split}
\end{equation}
Since $\MM_\ell \subseteq \TT_\ell \backslash \TT_\ell^{\rm bad} = \TT_\ell^{\rm good} \cup \TT_\ell^{\rm neither}$,
it follows from~\eqref{hypsing:eq:msv:step1}--\eqref{hypsing:eq:msv:step2} that
\begin{equation*}
 \theta \,\tau_\ell(T)^2
 \reff{eq:doerfler}\le \sum_{T \in \MM_\ell} \tau_\ell(T)^2
 \le \sum_{T \in \TT_\ell^{\rm good}} \tau_\ell(T)^2
 + \sum_{T \in \TT_\ell^{\rm neither}} \tau_\ell(T)^2
 \xrightarrow{\ell \to \infty} 0.
\end{equation*}
This concludes the proof.
\end{proof}

\bibliographystyle{alpha}
\bibliography{literature}

\end{document}